\documentclass[12pt,a4paper]{article}
\usepackage{theorem,enumerate}
\usepackage{amsmath,latexsym,amssymb,amsfonts}
\usepackage{eucal}
\usepackage{color}
\usepackage{comment}
\usepackage{eufrak}
\newcounter{Scounter}
\theorembodyfont{\normalfont\slshape}
\newtheorem{thm}{Theorem}[section]

\newtheorem{Thm}{Theorem}

\newtheorem{prop}[thm]{Proposition}
\newtheorem{lem}[thm]{Lemma}

\newtheorem{Q}[thm]{Question} 

\newtheorem{claim}{Claim}[section] 

\newtheorem{fact}[claim]{Fact}

\numberwithin{equation}{section}

\newcommand{\proof}{\medbreak\noindent\textit{Proof.}\quad}

\newcommand{\qed}{{$\quad\square$\vs{3.6}}}

\newcommand{\vs}[1]{\vspace*{#1 mm}}

\def\A{{ \mathcal{A}}}

\def\C{{ \mathcal{C}}}
\def\D{{ \mathcal{D}}}

\def\G{{ \mathcal{G}}}
\def\H{{ \mathcal{H}}}

\def\L{{ \mathcal{L}}}

\bfseries\normalfont

\title{Path-factors involving paths of order seven and nine}
\author{
Yoshimi Egawa$^1$ \and\
Michitaka Furuya$^1$\footnote{\texttt{e-mail:michitaka.furuya@gmail.com}} \vs{5}\\
$^1$\textsl{Department of Mathematical Information Science,} \\
\textsl{Tokyo University of Science,}\\
\textsl{1-3 Kagurazaka, Shinjuku-ku, Tokyo 162-8601, Japan }\\
}

\date{}

\begin{document}

\maketitle

\begin{abstract}
In this paper, we show the following two theorems (here $c_{i}(G-X)$ is the number of components $C$ of $G-X$ with $|V(C)|=i$):
(i)~If a graph $G$ satisfies $c_{1}(G-X)+\frac{1}{3}c_{3}(G-X)+\frac{1}{3}c_{5}(G-X)\leq \frac{2}{3}|X|$ for all $X\subseteq V(G)$, then $G$ has a $\{P_{2},P_{7}\}$-factor.
(ii)~If a graph $G$ satisfies $c_{1}(G-X)+c_{3}(G-X)+\frac{2}{3}c_{5}(G-X)+\frac{1}{3}c_{7}(G-X)\leq \frac{2}{3}|X|$ for all $X\subseteq V(G)$, then $G$ has a $\{P_{2},P_{9}\}$-factor.
\end{abstract}

\noindent
{\it Key words and phrases.}
path-factor, matching, hypomatchable graph.

\noindent
{\it AMS 2010 Mathematics Subject Classification.}
05C70.

%%%%%%%%%%%%%%%%%%%%%%%%%%%%%%%%%%%%%%%%%%%%%%%%%%%%%%%%%%%%%%%%%%%%%%%%%%%%%%%%%%%%%%%%%%%%%%%%%%%%%%%%%%%%%%%%%%%%%%%%
%%%%%%%%%%%%%%%%%%%%%%%%%%%%%%%%%%%%%%%%%%%%%%%%%%%%%%%%%%%%%%%%%%%%%%%%%%%%%%%%%%%%%%%%%%%%%%%%%%%%%%%%%%%%%%%%%%%%%%%%
%%%%%%%%%%%%%%%%%%%%%%%%%%%%%%%%%%%%%%%%%%%%%%%%%%%%%%%%%%%%%%%%%%%%%%%%%%%%%%%%%%%%%%%%%%%%%%%%%%%%%%%%%%%%%%%%%%%%%%%%
\section{Introduction}\label{sec1}
%%%%%%%%%%%%%%%%%%%%%%%%%%%%%%%%%%%%%%%%%%%%%%%%%%%%%%%%%%%%%%%%%%%%%%%%%%%%%%%%%%%%%%%%%%%%%%%%%%%%%%%%%%%%%%%%%%%%%%%%
%%%%%%%%%%%%%%%%%%%%%%%%%%%%%%%%%%%%%%%%%%%%%%%%%%%%%%%%%%%%%%%%%%%%%%%%%%%%%%%%%%%%%%%%%%%%%%%%%%%%%%%%%%%%%%%%%%%%%%%%
%%%%%%%%%%%%%%%%%%%%%%%%%%%%%%%%%%%%%%%%%%%%%%%%%%%%%%%%%%%%%%%%%%%%%%%%%%%%%%%%%%%%%%%%%%%%%%%%%%%%%%%%%%%%%%%%%%%%%%%%

In this paper, all graphs are finite and simple.
Let $G$ be a graph.
We let $V(G)$ and $E(G)$ denote the vertex set and the edge set of $G$, respectively.
For $u\in V(G)$, we let $N_{G}(u)$ and $d_{G}(u)$ denote the {\it neighborhood} and the {\it degree} of $u$, respectively.
For $U\subseteq V(G)$, we let $N_{G}(U)=(\bigcup _{u\in U}N_{G}(u))-U$.
For disjoint sets $X,Y\subseteq V(G)$, we let $E_{G}(X,Y)$ denote the set of edges of $G$ joining a vertex in $X$ and a vertex in $Y$.
For $X\subseteq V(G)$, we let $G[X]$ denote the subgraph of $G$ induced by $X$.
%For a graph $H$ and an integer $s\geq 2$, we let $sH$ denote the disjoint union of $s$ copies of $H$.
For two graphs $H_{1}$ and $H_{2}$, we let $H_{1}\cup H_{2}$ and $H_{1}+H_{2}$ denote the {\it union} and the {\it join} of $H_{1}$ and $H_{2}$, respectively.
For a graph $H$ and an integer $s\geq 2$, we let $sH$ denote the disjoint union of $s$ copies of $H$.
Let $K_{n}$ and $P_{n}$ denote the {\it complete graph} and the {\it path} of order $n$, respectively.
For terms and symbols not defined here, we refer the reader to \cite{D}.

Let again $G$ be a graph.
A subset $M$ of $E(G)$ is a {\it matching} if no two distinct edges in $M$ have a common endvertex.
If there is no fear of confusion, we often identify a matching $M$ of $G$ with the subgraph of $G$ induced by $M$.
A matching $M$ of $G$ is {\it perfect} if $V(M)=V(G)$.
For a set $\H$ of connected graphs, a spanning subgraph $F$ of $G$ is called an {\it $\H$-factor} if each component of $F$ is isomorphic to a graph in $\H$.
Note that a perfect matching can be regarded as a $\{P_{2}\}$-factor.
A {\it path-factor} of $G$ is a spanning subgraph whose components are paths of order at least $2$.
Since every path of order at least $2$ can be partitioned into paths of orders $2$ and $3$, a graph has a path-factor if and only if it has a $\{P_{2},P_{3}\}$-factor.
Akiyama, Avis and Era~\cite{AAE} gave a necessary and sufficient condition for the existence of a path-factor (here $i(G)$ denotes the number of isolated vertices of a graph $G$).

\begin{Thm}[Akiyama, Avis and Era~\cite{AAE}]%%%%%%%%%%%%%%%%%%%%%%%%%%%%%%%%%%%%%%%%%%%%%%%%%%%%%%%%%%%%%%%%%%%%%%%%%%%
\label{ThmA}
A graph $G$ has a $\{P_{2},P_{3}\}$-factor if and only if $i(G-X)\leq 2|X|$ for all $X\subseteq V(G)$.
\end{Thm}
%%%%%%%%%%%%%%%%%%%%%%%%%%%%%%%%%%%%%%%%%%%%%%%%%%%%%%%%%%%%%%%%%%%%%%%%%%%%%%%%%%%%%%%%%%%%%%%%%%%%%%%%%%%%%%%%%%%%%%%%

On the other hand, it follows from a result of Loebal and Poljak~\cite{LP} that for $k\geq 2$, the existence problem of a $\{P_{2},P_{2k+1}\}$-factor is {\bf NP}-complete.
However, in general, the fact that a problem is {\bf NP}-complete in terms of algorithm does not mean that one cannot obtain a theoretical result concerning the problem.
In this paper, we discuss sufficient conditions for the existence of a $\{P_{2},P_{2k+1}\}$-factor (for detailed historical background and motivations, we refer the reader to \cite{EF}).
%From this viewpoint, we focus on the existence of a $\{P_{2},P_{2k+1}\}$-factor.

In order to state our results, we need some more preparations.
For a graph $H$, we let $\C(H)$ be the set of components of $H$, and for $i\geq 1$, let $\C _{i}(H)=\{C\in \C(H)\mid |V(C)|=i\}$ and $c_{i}(H)=|\C _{i}(H)|$.
Note that $c_{1}(H)$ is the number of isolated vertices of $H$ (i.e., $c_{1}(H)=i(H)$).
%For $k\geq 1$, if a graph $G$ has a $\{P_{2},P_{2k+1}\}$-factor, then $\sum _{0\leq i\leq k-1}\frac{k-i}{k}c_{2i+1}(G-X)\leq \frac{k+1}{k}|X|$ for all $X\subseteq V(G)$ (see Section~\ref{sec2}).
For $k\geq 1$, if a graph $G$ has a $\{P_{2},P_{2k+1}\}$-factor, then $\sum _{0\leq i\leq k-1}(k-i)c_{2i+1}(G-X)\leq (k+1)|X|$ for all $X\subseteq V(G)$ (see Section~\ref{sec2}).
Thus if a condition concerning $c_{2i+1}(G-X)~(0\leq i\leq k-1)$ for $X\subseteq V(G)$ assures us the existence of a $\{P_{2},P_{2k+1}\}$-factor, then it will make a useful sufficient condition.

Recently, in \cite{EF}, the authors proved the following theorem, and showed that the bound $\frac{4}{3}|X|+\frac{1}{3}$ in the theorem is best possible.

\begin{Thm}[Egawa and Furuya~\cite{EF}]%%%%%%%%%%%%%%%%%%%%%%%%%%%%%%%%%%%%%%%%%%%%%%%%%%%%%%%%%%%%%%%%%%%%%%%%%%%
\label{ThmB}
Let $G$ be a graph.
If $c_{1}(G-X)+\frac{2}{3}c_{3}(G-X)\leq \frac{4}{3}|X|+\frac{1}{3}$ for all $X\subseteq V(G)$, then $G$ has a $\{P_{2},P_{5}\}$-factor.
\end{Thm}
%%%%%%%%%%%%%%%%%%%%%%%%%%%%%%%%%%%%%%%%%%%%%%%%%%%%%%%%%%%%%%%%%%%%%%%%%%%%%%%%%%%%%%%%%%%%%%%%%%%%%%%%%%%%%%%%%%%%%%%%

In \cite{EF}, the authors also constructed examples which show that for $k\geq 3$ with $k\equiv 0~(\mbox{mod }3)$, there exist infinitely many graphs $G$ having no $\{P_{2},P_{2k+1}\}$-factor such that $\sum _{0\leq i\leq k-1}c_{2i+1}(G-X)\leq \frac{4k+6}{8k+3}|X|+\frac{2k+3}{8k+3}$ for all $X\subseteq V(G)$, and proposed a conjecture that, for an integer $k\geq 3$ and a graph $G$, if $\sum _{0\leq i\leq k-1}c_{2i+1}(G-X)\leq \frac{4k+6}{8k+3}|X|$ for all $X\subseteq V(G)$, then $G$ has a $\{P_{2},P_{2k+1}\}$-factor.

In this paper, we settle the above conjecture for the case where $k\in \{3,4\}$ as follows (note that Theorem~\ref{mainthm2} implies that the coefficient $\frac{4k+6}{8k+3}$ of $|X|$ in the conjecture is not best possible for $k=4$).

\begin{thm}%%%%%%%%%%%%%%%%%%%%%%%%%%%%%%%%%%%%%%%%%%%%%%%%%%%%%%%%%%%%%%%%%%%%%%%%%%%%%%%%%%%%%%%%%%%%%%%%%%%%%%%%%%%%%
\label{mainthm1}
Let $G$ be a graph.
If $c_{1}(G-X)+\frac{1}{3}c_{3}(G-X)+\frac{1}{3}c_{5}(G-X)\leq \frac{2}{3}|X|$ for all $X\subseteq V(G)$, then $G$ has a $\{P_{2},P_{7}\}$-factor.
\end{thm}
%%%%%%%%%%%%%%%%%%%%%%%%%%%%%%%%%%%%%%%%%%%%%%%%%%%%%%%%%%%%%%%%%%%%%%%%%%%%%%%%%%%%%%%%%%%%%%%%%%%%%%%%%%%%%%%%%%%%%%%%

\begin{thm}%%%%%%%%%%%%%%%%%%%%%%%%%%%%%%%%%%%%%%%%%%%%%%%%%%%%%%%%%%%%%%%%%%%%%%%%%%%%%%%%%%%%%%%%%%%%%%%%%%%%%%%%%%%%%
\label{mainthm2}
Let $G$ be a graph.
If $c_{1}(G-X)+c_{3}(G-X)+\frac{2}{3}c_{5}(G-X)+\frac{1}{3}c_{7}(G-X)\leq \frac{2}{3}|X|$ for all $X\subseteq V(G)$, then $G$ has a $\{P_{2},P_{9}\}$-factor.
\end{thm}
%%%%%%%%%%%%%%%%%%%%%%%%%%%%%%%%%%%%%%%%%%%%%%%%%%%%%%%%%%%%%%%%%%%%%%%%%%%%%%%%%%%%%%%%%%%%%%%%%%%%%%%%%%%%%%%%%%%%%%%%

We prove Theorems~\ref{mainthm1} and \ref{mainthm2} in Sections~\ref{sec3}--\ref{sec5}.
We remark that hypomatchable graphs play an important role in the proof, through $P_{7}$ and $P_{9}$ are not hypomatchable (see Section~\ref{sec4} for the definition of a hypomatchable graph).
In Section~\ref{sec6}, we discuss the sharpness of coefficients in Theorems~\ref{mainthm1} and \ref{mainthm2}.

In our proof of Theorems~\ref{mainthm1} and \ref{mainthm2}, we make use of the following fact.

\begin{fact}%%%%%%%%%%%%%%%%%%%%%%%%%%%%%%%%%%%%%%%%%%%%%%%%%%%%%%%%%%%%%%%%%%%%%%%%%%%%%%%%%%%%%%%%%%%%%%%%%%%%%%%%%%%%
\label{fact1}
Let $k\geq 2$ be an integer, and let $G$ be a graph.
Then $G$ has a $\{P_{2},P_{2k+1}\}$-factor if and only if $G$ has a path-factor $F$ such that $\C_{2i+1}(F)=\emptyset $ for every $i~(1\leq i\leq k-1)$.
\end{fact}
%%%%%%%%%%%%%%%%%%%%%%%%%%%%%%%%%%%%%%%%%%%%%%%%%%%%%%%%%%%%%%%%%%%%%%%%%%%%%%%%%%%%%%%%%%%%%%%%%%%%%%%%%%%%%%%%%%%%%%%%

%%%%%%%%%%%%%%%%%%%%%%%%%%%%%%%%%%%%%%%%%%%%%%%%%%%%%%%%%%%%%%%%%%%%%%%%%%%%%%%%%%%%%%%%%%%%%%%%%%%%%%%%%%%%%%%%%%%%%%%%
%%%%%%%%%%%%%%%%%%%%%%%%%%%%%%%%%%%%%%%%%%%%%%%%%%%%%%%%%%%%%%%%%%%%%%%%%%%%%%%%%%%%%%%%%%%%%%%%%%%%%%%%%%%%%%%%%%%%%%%%
%%%%%%%%%%%%%%%%%%%%%%%%%%%%%%%%%%%%%%%%%%%%%%%%%%%%%%%%%%%%%%%%%%%%%%%%%%%%%%%%%%%%%%%%%%%%%%%%%%%%%%%%%%%%%%%%%%%%%%%%
\section{A necessary condition for $\{P_{2},P_{2k+1}\}$-factor}\label{sec2}
%%%%%%%%%%%%%%%%%%%%%%%%%%%%%%%%%%%%%%%%%%%%%%%%%%%%%%%%%%%%%%%%%%%%%%%%%%%%%%%%%%%%%%%%%%%%%%%%%%%%%%%%%%%%%%%%%%%%%%%%
%%%%%%%%%%%%%%%%%%%%%%%%%%%%%%%%%%%%%%%%%%%%%%%%%%%%%%%%%%%%%%%%%%%%%%%%%%%%%%%%%%%%%%%%%%%%%%%%%%%%%%%%%%%%%%%%%%%%%%%%
%%%%%%%%%%%%%%%%%%%%%%%%%%%%%%%%%%%%%%%%%%%%%%%%%%%%%%%%%%%%%%%%%%%%%%%%%%%%%%%%%%%%%%%%%%%%%%%%%%%%%%%%%%%%%%%%%%%%%%%%

In this section, we give a necessary condition for the existence of a $\{P_{2},P_{2k+1}\}$-factor in terms of invariants $c_{2i+1}~(0\leq i\leq k-1)$.
We show the following proposition.

\begin{prop}%%%%%%%%%%%%%%%%%%%%%%%%%%%%%%%%%%%%%%%%%%%%%%%%%%%%%%%%%%%%%%%%%%%%%%%%%%%%%%%%%%%%%%%%%%%%%%%%%%%%%%%%%%%%
\label{prop2.1}
For an integer $k\geq 1$, if a graph $G$ has a $\{P_{2},P_{2k+1}\}$-factor, then $\sum _{0\leq i\leq k-1}(k-i)c_{2i+1}(G-X)\leq (k+1)|X|$ for all $X\subseteq V(G)$.
\end{prop}
%%%%%%%%%%%%%%%%%%%%%%%%%%%%%%%%%%%%%%%%%%%%%%%%%%%%%%%%%%%%%%%%%%%%%%%%%%%%%%%%%%%%%%%%%%%%%%%%%%%%%%%%%%%%%%%%%%%%%%%%
\proof
Let $F$ be a $\{P_{2},P_{2k+1}\}$-factor of $G$, and let $X\subseteq V(G)$.
Observe that
$$
\sum _{0\leq i\leq k-1}(k-i)c_{2i+1}(G-X)=\sum _{C\in \bigcup _{0\leq i\leq k-1}\C_{2i+1}(G-X)}\left(k+\frac{1}{2}-\frac{|V(C)|}{2}\right).
$$
%and
%$$
%\sum _{0\leq i\leq k-1}(k-i)c_{2i+1}(F-X)=\sum _{H\in \bigcup _{0\leq i\leq k-1}\C_{2i+1}(F-X)}\left(k-\frac{|V(H)|-1}{2}\right).
%$$
With this observation in mind, we first prove the following claim.

\begin{claim}%%%%%%%%%%%%%%%%%%%%%%%%%%%%%%%%%%%%%%%%%%%%%%%%%%%%%%%%%%%%%%%%%%%%%%%%%%%%%%%%%%%%%%%%%%%%%%%%%%%%%%%%%%%
\label{cl2.1.1}
Let $P\in \C(F)$.
Then $\sum _{H\in \bigcup _{0\leq i\leq k-1}\C_{2i+1}(P-Y)}(k+\frac{1}{2}-\frac{|V(H)|}{2})\leq (k+1)|Y|$ for all $Y\subseteq V(P)$.
\end{claim}
%%%%%%%%%%%%%%%%%%%%%%%%%%%%%%%%%%%%%%%%%%%%%%%%%%%%%%%%%%%%%%%%%%%%%%%%%%%%%%%%%%%%%%%%%%%%%%%%%%%%%%%%%%%%%%%%%%%%%%%%
\proof
We proceed by induction on $|Y|$.
If $Y=\emptyset $, the desired inequality clearly holds.
Thus let $Y\not=\emptyset $, and assume that the desired inequality holds for subsets of $V(P)$ with cardinality $|Y|-1$.
Take $x\in Y$, and set $Y'=Y-\{x\}$.
Then $\sum _{H\in \bigcup _{0\leq i\leq k-1}\C_{2i+1}(P-Y')}(k+\frac{1}{2}-\frac{|V(H)|}{2})\leq (k+1)|Y'|$.
Let $H_{0}$ be the component of $P-Y'$ containing $x$, and let $H_{1}$ and $H_{2}$ denote the two segments of $H_{0}$ obtained by deleting $x$ from $H_{0}$.
Note that $H_{1}$ or $H_{2}$ (or both) may be empty.
If $H_{0}$ has even order, then precisely one of $H_{1}$ and $H_{2}$, say $H_{1}$, has odd order, and hence
\begin{align*}
\sum _{H\in \bigcup _{0\leq i\leq k-1}\C_{2i+1}(P-Y)}&\left(k+\frac{1}{2}-\frac{|V(H)|}{2}\right)\\
&= \sum _{H\in \bigcup _{0\leq i\leq k-1}\C_{2i+1}(P-Y')}\left(k+\frac{1}{2}-\frac{|V(H)|}{2}\right)+\left(k+\frac{1}{2}-\frac{|V(H_{1})|}{2}\right)\\
&\leq (k+1)|Y'|+k\\
&<(k+1)|Y|.
\end{align*}
Thus we may assume that $H_{0}$ has odd order.
Note that $-(k+\frac{1}{2}-\frac{|V(H_{0})|}{2})+(k+\frac{1}{2}-\frac{|V(H_{1})|}{2})+(k+\frac{1}{2}-\frac{|V(H_{2})|}{2})=k+\frac{1}{2}+\frac{|V(H_{0})|-|V(H_{1})|-|V(H_{2})|}{2}=k+1$.
Consequently
\begin{align*}
\sum _{H\in \bigcup _{0\leq i\leq k-1}\C_{2i+1}(P-Y)}&\left(k+\frac{1}{2}-\frac{|V(H)|}{2}\right)\\
&\leq \sum _{H\in \bigcup _{0\leq i\leq k-1}\C_{2i+1}(P-Y')}\left(k+\frac{1}{2}-\frac{|V(H)|}{2}\right)-\left(k+\frac{1}{2}-\frac{|V(H_{0})|}{2}\right)\\
&\quad +\left(k+\frac{1}{2}-\frac{|V(H_{1})|}{2}\right)+\left(k+\frac{1}{2}-\frac{|V(H_{2})|}{2}\right)\\
&\leq (k+1)|Y'|+(k+1)\\
&=(k+1)|Y|,
\end{align*}
as desired (note that this argument works even if $Y'=\emptyset $ and $H_{0}=P$).
\qed

Let $C\in \bigcup _{0\leq i\leq k-1}\C_{2i+1}(G-X)$.
Since $|V(C)|$ is odd, $F[V(C)]$ has a component $H_{C}$ of odd order.
We have $|V(H_{C})|\leq |V(C)|$ and $H_{C}\in \bigcup _{0\leq i\leq k-1}\C_{2i+1}(F-X)$.
Now let $\H=\{H_{C}\mid C\in \bigcup _{0\leq i\leq k-1}\C_{2i+1}(G-X)\}$.
Clearly we have $H_{C}\not=H_{C'}$ for any $C,C'\in \bigcup _{0\leq i\leq k-1}\C_{2i+1}(G-X)$ with $C\not=C'$.
Consequently
\begin{align*}
\sum _{0\leq i\leq k-1}(k-i)c_{2i+1}(G-X) &= \sum _{C\in \bigcup _{0\leq i\leq k-1}\C_{2i+1}(G-X)}\left(k+\frac{1}{2}-\frac{|V(C)|}{2}\right)\\
&\leq \sum _{C\in \bigcup _{0\leq i\leq k-1}\C_{2i+1}(G-X)}\left(k+\frac{1}{2}-\frac{|V(H_{C})|}{2}\right)\\
&=\sum _{H\in \H}\left(k+\frac{1}{2}-\frac{|V(H)|}{2}\right)\\
&\leq \sum _{H\in \bigcup _{0\leq i\leq k-1}\C_{2i+1}(F-X)}\left(k+\frac{1}{2}-\frac{|V(H)|}{2}\right)\\
&=\sum _{P\in \C(F)}\left(\sum _{H\in \bigcup _{0\leq i\leq k-1}\C_{2i+1}(P-X)}\left(k+\frac{1}{2}-\frac{|V(H)|}{2}\right)\right).
\end{align*}
Therefore it follows from Claim~\ref{cl2.1.1} that
\begin{align*}
\sum _{0\leq i\leq k-1}(k-i)c_{2i+1}(G-X)&\leq \sum _{P\in \C(F)}(k+1)|V(P)\cap X|\\
&= (k+1)|X|,
\end{align*}
as desired.
\qed

%%%%%%%%%%%%%%%%%%%%%%%%%%%%%%%%%%%%%%%%%%%%%%%%%%%%%%%%%%%%%%%%%%%%%%%%%%%%%%%%%%%%%%%%%%%%%%%%%%%%%%%%%%%%%%%%%%%%%%%%
%%%%%%%%%%%%%%%%%%%%%%%%%%%%%%%%%%%%%%%%%%%%%%%%%%%%%%%%%%%%%%%%%%%%%%%%%%%%%%%%%%%%%%%%%%%%%%%%%%%%%%%%%%%%%%%%%%%%%%%%
%%%%%%%%%%%%%%%%%%%%%%%%%%%%%%%%%%%%%%%%%%%%%%%%%%%%%%%%%%%%%%%%%%%%%%%%%%%%%%%%%%%%%%%%%%%%%%%%%%%%%%%%%%%%%%%%%%%%%%%%
\section{Linear forests in bipartite graphs}\label{sec3}
%%%%%%%%%%%%%%%%%%%%%%%%%%%%%%%%%%%%%%%%%%%%%%%%%%%%%%%%%%%%%%%%%%%%%%%%%%%%%%%%%%%%%%%%%%%%%%%%%%%%%%%%%%%%%%%%%%%%%%%%
%%%%%%%%%%%%%%%%%%%%%%%%%%%%%%%%%%%%%%%%%%%%%%%%%%%%%%%%%%%%%%%%%%%%%%%%%%%%%%%%%%%%%%%%%%%%%%%%%%%%%%%%%%%%%%%%%%%%%%%%
%%%%%%%%%%%%%%%%%%%%%%%%%%%%%%%%%%%%%%%%%%%%%%%%%%%%%%%%%%%%%%%%%%%%%%%%%%%%%%%%%%%%%%%%%%%%%%%%%%%%%%%%%%%%%%%%%%%%%%%%

In this this section, we show the following proposition, which plays a key role in the proof of our main theorems.

\begin{prop}%%%%%%%%%%%%%%%%%%%%%%%%%%%%%%%%%%%%%%%%%%%%%%%%%%%%%%%%%%%%%%%%%%%%%%%%%%%%%%%%%%%%%%%%%%%%%%%%%%%%%%%%%%%%
\label{prop3.1}
Let $S$ and $T$ be disjoint sets, and let $T_{1}$ and $T_{2}$ be disjoint subsets of $T$.
Let $G$ be a bipartite graph with bipartition $(S,T)$, and let $L\subseteq E(G)$.
Suppose that
\begin{enumerate}[{\upshape(i)}]
\item
$|N_{G}(X)|\geq |X|$ for every $X\subseteq S$, and
\item
$|N_{G-L}(Y)|\geq |Y\cap T_{1}|+\frac{1}{2}|Y\cap T_{2}|$ for every $Y\subseteq T_{1}\cup T_{2}$.
\end{enumerate}
Then $G$ has a subgraph $F$ with $V(F)\supseteq S\cup T_{1}\cup T_{2}$ such that each $A\in \C(F)$ is a path satisfying one of the following two conditions:
\begin{enumerate}[{\upshape(I)}]
\item
$|V(A)|=2$; or
\item
$E(A)\subseteq E(G)-L$, $V(A)\cap T\subseteq T_{1}\cup T_{2}$, $|V(A)\cap T_{2}|=2$ and the two vertices in $V(A)\cap T_{2}$ are the endvertices of $A$.
%$E(A)\subseteq E(G)-L$, $V(A)\cap T\subseteq T_{1}\cap T_{2}$, $|V(A)\cap T_{2}|=2$ and the two vertices in $V(A)\cap T_{2}$ are the endvertices of $A$.
\end{enumerate}
\end{prop}
%%%%%%%%%%%%%%%%%%%%%%%%%%%%%%%%%%%%%%%%%%%%%%%%%%%%%%%%%%%%%%%%%%%%%%%%%%%%%%%%%%%%%%%%%%%%%%%%%%%%%%%%%%%%%%%%%%%%%%%%

As a preparation for the proof of Proposition~\ref{prop3.1}, we first show the following lemma.

\begin{lem}%%%%%%%%%%%%%%%%%%%%%%%%%%%%%%%%%%%%%%%%%%%%%%%%%%%%%%%%%%%%%%%%%%%%%%%%%%%%%%%%%%%%%%%%%%%%%%%%%%%%%%%%%%%%%
\label{lem3.1.1}
Let $S$ and $T$ be disjoint sets, and let $T_{1}$ and $T_{2}$ be disjoint subsets of $T$ such that $T_{1}\cup T_{2}=T$.
Let $H$ be a bipartite graph with bipartition $(S,T)$, and suppose that $|N_{H}(Y)|\geq |Y\cap T_{1}|+\frac{1}{2}|Y\cap T_{2}|$ for every $Y\subseteq T$.
Then $H$ has a subgraph $F$ with $V(F)\supseteq T_{1}\cup T_{2}$ such that each $A\in \C(F)$ is a path satisfying one of the following two conditions:
\begin{enumerate}[{\upshape (I')}]
\item
$|V(A)|=2$; or
\item
$|V(A)\cap T_{2}|=2$ and the two vertices in $V(A)\cap T_{2}$ are the endvertices of $A$.
\end{enumerate}
\end{lem}
%%%%%%%%%%%%%%%%%%%%%%%%%%%%%%%%%%%%%%%%%%%%%%%%%%%%%%%%%%%%%%%%%%%%%%%%%%%%%%%%%%%%%%%%%%%%%%%%%%%%%%%%%%%%%%%%%%%%%%%%
\proof
By the assumption of the lemma, $|N_{H}(Y)|\geq |Y\cap T_{1}|+\frac{1}{2}|Y\cap T_{2}|=|Y|$ for every $Y\subseteq T_{1}$.
Hence by Hall's marriage theorem, there exists a matching $F$ of $H$ such that $V(F)\cap T=T_{1}$.
In particular, $H$ has a subgraph $F$ with $V(F)\supseteq T_{1}$ such that each $A\in \C(F)$ is a path satisfying (I') or (II').
Choose such a subgraph $F$ so that $|(S\cup T_{2})-V(F)|$ is as small as possible.

It suffices to show that $T_{2}-V(F)=\emptyset $.
By way of contradiction, suppose that $T_{2}-V(F)\not=\emptyset $.
Now we define the set $\A$ of paths of $H$ as follows:
Let $\A_{0}$ be the set of paths of $H$ consisting of one vertex in $T_{2}-V(F)$.
For each $i\geq 1$, let $\A_{i}$ be the set of components $A$ of $F$ with $A\not\in \bigcup _{0\leq j\leq i-1}\A_{j}$ and $E_{H}(V(A)\cap S,\bigcup _{A'\in \A_{i-1}}(V(A')\cap T))\not=\emptyset $.
Let $\A=\bigcup _{i\geq 0}\A_{i}$.

\begin{figure}
\begin{center}
%WinTpicVersion4.28b
{\unitlength 0.1in
\begin{picture}( 46.0000,  9.0000)(  2.0000,-11.3500)
% CIRCLE 2 0 0 0 Black Black
% 4 595 600 595 650 595 650 595 650
% 
\special{sh 1.000}%
\special{ia 596 600 50 50  0.0000000  6.2831853}%
\special{pn 8}%
\special{ar 596 600 50 50  0.0000000  6.2831853}%
% CIRCLE 2 0 0 0 Black Black
% 4 695 1000 695 1050 695 1050 695 1050
% 
\special{sh 1.000}%
\special{ia 696 1000 50 50  0.0000000  6.2831853}%
\special{pn 8}%
\special{ar 696 1000 50 50  0.0000000  6.2831853}%
% CIRCLE 2 0 0 0 Black Black
% 4 795 600 795 650 795 650 795 650
% 
\special{sh 1.000}%
\special{ia 796 600 50 50  0.0000000  6.2831853}%
\special{pn 8}%
\special{ar 796 600 50 50  0.0000000  6.2831853}%
% CIRCLE 2 0 0 0 Black Black
% 4 995 600 995 650 995 650 995 650
% 
\special{sh 1.000}%
\special{ia 996 600 50 50  0.0000000  6.2831853}%
\special{pn 8}%
\special{ar 996 600 50 50  0.0000000  6.2831853}%
% CIRCLE 2 0 0 0 Black Black
% 4 1195 600 1195 650 1195 650 1195 650
% 
\special{sh 1.000}%
\special{ia 1196 600 50 50  0.0000000  6.2831853}%
\special{pn 8}%
\special{ar 1196 600 50 50  0.0000000  6.2831853}%
% CIRCLE 2 0 0 0 Black Black
% 4 895 1000 895 1050 895 1050 895 1050
% 
\special{sh 1.000}%
\special{ia 896 1000 50 50  0.0000000  6.2831853}%
\special{pn 8}%
\special{ar 896 1000 50 50  0.0000000  6.2831853}%
% CIRCLE 2 0 0 0 Black Black
% 4 1095 1000 1095 1050 1095 1050 1095 1050
% 
\special{sh 1.000}%
\special{ia 1096 1000 50 50  0.0000000  6.2831853}%
\special{pn 8}%
\special{ar 1096 1000 50 50  0.0000000  6.2831853}%
% CIRCLE 2 0 0 0 Black Black
% 4 1395 600 1395 650 1395 650 1395 650
% 
\special{sh 1.000}%
\special{ia 1396 600 50 50  0.0000000  6.2831853}%
\special{pn 8}%
\special{ar 1396 600 50 50  0.0000000  6.2831853}%
% CIRCLE 2 0 0 0 Black Black
% 4 1295 1000 1295 1050 1295 1050 1295 1050
% 
\special{sh 1.000}%
\special{ia 1296 1000 50 50  0.0000000  6.2831853}%
\special{pn 8}%
\special{ar 1296 1000 50 50  0.0000000  6.2831853}%
% LINE 1 0 3 0 Black Black
% 2 595 600 695 1000
% 
\special{pn 13}%
\special{pa 596 600}%
\special{pa 696 1000}%
\special{fp}%
% LINE 1 0 3 0 Black Black
% 2 795 600 895 1000
% 
\special{pn 13}%
\special{pa 796 600}%
\special{pa 896 1000}%
\special{fp}%
% LINE 1 0 3 0 Black Black
% 2 995 600 1095 1000
% 
\special{pn 13}%
\special{pa 996 600}%
\special{pa 1096 1000}%
\special{fp}%
% LINE 1 0 3 0 Black Black
% 2 1195 600 1295 1000
% 
\special{pn 13}%
\special{pa 1196 600}%
\special{pa 1296 1000}%
\special{fp}%
% LINE 1 0 3 0 Black Black
% 2 1395 600 1295 1000
% 
\special{pn 13}%
\special{pa 1396 600}%
\special{pa 1296 1000}%
\special{fp}%
% LINE 1 0 3 0 Black Black
% 2 1195 600 1095 1000
% 
\special{pn 13}%
\special{pa 1196 600}%
\special{pa 1096 1000}%
\special{fp}%
% LINE 1 0 3 0 Black Black
% 2 995 600 895 1000
% 
\special{pn 13}%
\special{pa 996 600}%
\special{pa 896 1000}%
\special{fp}%
% LINE 1 0 3 0 Black Black
% 2 795 600 695 1000
% 
\special{pn 13}%
\special{pa 796 600}%
\special{pa 696 1000}%
\special{fp}%
% CIRCLE 2 0 0 0 Black Black
% 4 1490 1000 1490 1050 1490 1050 1490 1050
% 
\special{sh 1.000}%
\special{ia 1490 1000 50 50  0.0000000  6.2831853}%
\special{pn 8}%
\special{ar 1490 1000 50 50  0.0000000  6.2831853}%
% CIRCLE 2 0 0 0 Black Black
% 4 1590 600 1590 650 1590 650 1590 650
% 
\special{sh 1.000}%
\special{ia 1590 600 50 50  0.0000000  6.2831853}%
\special{pn 8}%
\special{ar 1590 600 50 50  0.0000000  6.2831853}%
% CIRCLE 2 0 0 0 Black Black
% 4 1790 600 1790 650 1790 650 1790 650
% 
\special{sh 1.000}%
\special{ia 1790 600 50 50  0.0000000  6.2831853}%
\special{pn 8}%
\special{ar 1790 600 50 50  0.0000000  6.2831853}%
% CIRCLE 2 0 0 0 Black Black
% 4 1690 1000 1690 1050 1690 1050 1690 1050
% 
\special{sh 1.000}%
\special{ia 1690 1000 50 50  0.0000000  6.2831853}%
\special{pn 8}%
\special{ar 1690 1000 50 50  0.0000000  6.2831853}%
% LINE 1 0 3 0 Black Black
% 2 1390 600 1490 1000
% 
\special{pn 13}%
\special{pa 1390 600}%
\special{pa 1490 1000}%
\special{fp}%
% LINE 1 0 3 0 Black Black
% 2 1590 600 1690 1000
% 
\special{pn 13}%
\special{pa 1590 600}%
\special{pa 1690 1000}%
\special{fp}%
% LINE 1 0 3 0 Black Black
% 2 1790 600 1690 1000
% 
\special{pn 13}%
\special{pa 1790 600}%
\special{pa 1690 1000}%
\special{fp}%
% LINE 1 0 3 0 Black Black
% 2 1590 600 1490 1000
% 
\special{pn 13}%
\special{pa 1590 600}%
\special{pa 1490 1000}%
\special{fp}%
% SPLINE 2 2 3 0 Black Black
% 4 1000 600 800 400 200 510 200 510
% 
\special{pn 8}%
\special{pn 8}%
\special{pa 1000 600}%
\special{pa 995 594}%
\special{fp}%
\special{pa 972 564}%
\special{pa 967 558}%
\special{fp}%
\special{pa 944 529}%
\special{pa 938 523}%
\special{fp}%
\special{pa 913 495}%
\special{pa 908 489}%
\special{fp}%
\special{pa 882 462}%
\special{pa 876 456}%
\special{fp}%
\special{pa 848 432}%
\special{pa 842 427}%
\special{fp}%
\special{pa 811 406}%
\special{pa 804 402}%
\special{fp}%
\special{pa 771 386}%
\special{pa 763 384}%
\special{fp}%
\special{pa 727 374}%
\special{pa 719 373}%
\special{fp}%
\special{pa 682 368}%
\special{pa 674 368}%
\special{fp}%
\special{pa 637 368}%
\special{pa 629 368}%
\special{fp}%
\special{pa 592 373}%
\special{pa 584 374}%
\special{fp}%
\special{pa 547 381}%
\special{pa 539 383}%
\special{fp}%
\special{pa 503 392}%
\special{pa 495 394}%
\special{fp}%
\special{pa 459 405}%
\special{pa 452 407}%
\special{fp}%
\special{pa 416 419}%
\special{pa 409 422}%
\special{fp}%
\special{pa 374 435}%
\special{pa 367 439}%
\special{fp}%
\special{pa 332 453}%
\special{pa 325 456}%
\special{fp}%
\special{pa 291 471}%
\special{pa 283 474}%
\special{fp}%
\special{pa 249 489}%
\special{pa 242 492}%
\special{fp}%
\special{pa 207 507}%
\special{pa 200 510}%
\special{fp}%
% SPLINE 2 2 3 0 Black Black
% 6 1500 1000 1500 600 1700 440 2100 440 2200 510 2200 510
% 
\special{pn 8}%
\special{pn 8}%
\special{pa 1500 1000}%
\special{pa 1498 992}%
\special{fp}%
\special{pa 1491 955}%
\special{pa 1489 947}%
\special{fp}%
\special{pa 1482 910}%
\special{pa 1481 902}%
\special{fp}%
\special{pa 1475 864}%
\special{pa 1474 856}%
\special{fp}%
\special{pa 1469 819}%
\special{pa 1468 811}%
\special{fp}%
\special{pa 1466 773}%
\special{pa 1466 765}%
\special{fp}%
\special{pa 1466 727}%
\special{pa 1467 719}%
\special{fp}%
\special{pa 1472 682}%
\special{pa 1473 673}%
\special{fp}%
\special{pa 1484 637}%
\special{pa 1487 630}%
\special{fp}%
\special{pa 1504 596}%
\special{pa 1508 589}%
\special{fp}%
\special{pa 1530 558}%
\special{pa 1536 552}%
\special{fp}%
\special{pa 1562 526}%
\special{pa 1568 520}%
\special{fp}%
\special{pa 1598 497}%
\special{pa 1605 492}%
\special{fp}%
\special{pa 1637 472}%
\special{pa 1644 468}%
\special{fp}%
\special{pa 1678 451}%
\special{pa 1685 448}%
\special{fp}%
\special{pa 1720 433}%
\special{pa 1727 430}%
\special{fp}%
\special{pa 1763 418}%
\special{pa 1771 417}%
\special{fp}%
\special{pa 1808 407}%
\special{pa 1816 406}%
\special{fp}%
\special{pa 1853 400}%
\special{pa 1861 399}%
\special{fp}%
\special{pa 1899 397}%
\special{pa 1907 396}%
\special{fp}%
\special{pa 1945 398}%
\special{pa 1953 399}%
\special{fp}%
\special{pa 1991 404}%
\special{pa 1999 405}%
\special{fp}%
\special{pa 2035 414}%
\special{pa 2043 416}%
\special{fp}%
\special{pa 2079 429}%
\special{pa 2086 433}%
\special{fp}%
\special{pa 2119 451}%
\special{pa 2126 455}%
\special{fp}%
\special{pa 2157 477}%
\special{pa 2164 482}%
\special{fp}%
\special{pa 2194 505}%
\special{pa 2200 510}%
\special{fp}%
% STR 2 0 3 0 Black Black
% 4 1200 1100 1200 1200 5 0 0 0
% $A_{j}$
\put(12.0000,-12.0000){\makebox(0,0){$A_{j}$}}%
% VECTOR 0 0 3 0 Black Black
% 2 2300 800 2700 800
% 
\special{pn 20}%
\special{pa 2300 800}%
\special{pa 2700 800}%
\special{fp}%
\special{sh 1}%
\special{pa 2700 800}%
\special{pa 2634 780}%
\special{pa 2648 800}%
\special{pa 2634 820}%
\special{pa 2700 800}%
\special{fp}%
% CIRCLE 2 0 0 0 Black Black
% 4 3195 593 3195 643 3195 643 3195 643
% 
\special{sh 1.000}%
\special{ia 3196 594 50 50  0.0000000  6.2831853}%
\special{pn 8}%
\special{ar 3196 594 50 50  0.0000000  6.2831853}%
% CIRCLE 2 0 0 0 Black Black
% 4 3295 993 3295 1043 3295 1043 3295 1043
% 
\special{sh 1.000}%
\special{ia 3296 994 50 50  0.0000000  6.2831853}%
\special{pn 8}%
\special{ar 3296 994 50 50  0.0000000  6.2831853}%
% CIRCLE 2 0 0 0 Black Black
% 4 3395 593 3395 643 3395 643 3395 643
% 
\special{sh 1.000}%
\special{ia 3396 594 50 50  0.0000000  6.2831853}%
\special{pn 8}%
\special{ar 3396 594 50 50  0.0000000  6.2831853}%
% CIRCLE 2 0 0 0 Black Black
% 4 3595 593 3595 643 3595 643 3595 643
% 
\special{sh 1.000}%
\special{ia 3596 594 50 50  0.0000000  6.2831853}%
\special{pn 8}%
\special{ar 3596 594 50 50  0.0000000  6.2831853}%
% CIRCLE 2 0 0 0 Black Black
% 4 3795 593 3795 643 3795 643 3795 643
% 
\special{sh 1.000}%
\special{ia 3796 594 50 50  0.0000000  6.2831853}%
\special{pn 8}%
\special{ar 3796 594 50 50  0.0000000  6.2831853}%
% CIRCLE 2 0 0 0 Black Black
% 4 3495 993 3495 1043 3495 1043 3495 1043
% 
\special{sh 1.000}%
\special{ia 3496 994 50 50  0.0000000  6.2831853}%
\special{pn 8}%
\special{ar 3496 994 50 50  0.0000000  6.2831853}%
% CIRCLE 2 0 0 0 Black Black
% 4 3695 993 3695 1043 3695 1043 3695 1043
% 
\special{sh 1.000}%
\special{ia 3696 994 50 50  0.0000000  6.2831853}%
\special{pn 8}%
\special{ar 3696 994 50 50  0.0000000  6.2831853}%
% CIRCLE 2 0 0 0 Black Black
% 4 3995 593 3995 643 3995 643 3995 643
% 
\special{sh 1.000}%
\special{ia 3996 594 50 50  0.0000000  6.2831853}%
\special{pn 8}%
\special{ar 3996 594 50 50  0.0000000  6.2831853}%
% CIRCLE 2 0 0 0 Black Black
% 4 3895 993 3895 1043 3895 1043 3895 1043
% 
\special{sh 1.000}%
\special{ia 3896 994 50 50  0.0000000  6.2831853}%
\special{pn 8}%
\special{ar 3896 994 50 50  0.0000000  6.2831853}%
% LINE 0 0 3 0 Black Black
% 2 3195 593 3295 993
% 
\special{pn 20}%
\special{pa 3196 594}%
\special{pa 3296 994}%
\special{fp}%
% LINE 0 0 3 0 Black Black
% 2 3395 593 3495 993
% 
\special{pn 20}%
\special{pa 3396 594}%
\special{pa 3496 994}%
\special{fp}%
% LINE 2 2 3 0 Black Black
% 2 3595 593 3695 993
% 
\special{pn 8}%
\special{pa 3596 594}%
\special{pa 3696 994}%
\special{dt 0.045}%
% LINE 2 2 3 0 Black Black
% 2 3795 593 3895 993
% 
\special{pn 8}%
\special{pa 3796 594}%
\special{pa 3896 994}%
\special{dt 0.045}%
% LINE 2 2 3 0 Black Black
% 2 3995 593 3895 993
% 
\special{pn 8}%
\special{pa 3996 594}%
\special{pa 3896 994}%
\special{dt 0.045}%
% LINE 2 2 3 0 Black Black
% 2 3795 593 3695 993
% 
\special{pn 8}%
\special{pa 3796 594}%
\special{pa 3696 994}%
\special{dt 0.045}%
% LINE 0 0 3 0 Black Black
% 2 3595 593 3495 993
% 
\special{pn 20}%
\special{pa 3596 594}%
\special{pa 3496 994}%
\special{fp}%
% LINE 0 0 3 0 Black Black
% 2 3395 593 3295 993
% 
\special{pn 20}%
\special{pa 3396 594}%
\special{pa 3296 994}%
\special{fp}%
% CIRCLE 2 0 0 0 Black Black
% 4 4090 993 4090 1043 4090 1043 4090 1043
% 
\special{sh 1.000}%
\special{ia 4090 994 50 50  0.0000000  6.2831853}%
\special{pn 8}%
\special{ar 4090 994 50 50  0.0000000  6.2831853}%
% CIRCLE 2 0 0 0 Black Black
% 4 4190 593 4190 643 4190 643 4190 643
% 
\special{sh 1.000}%
\special{ia 4190 594 50 50  0.0000000  6.2831853}%
\special{pn 8}%
\special{ar 4190 594 50 50  0.0000000  6.2831853}%
% CIRCLE 2 0 0 0 Black Black
% 4 4390 593 4390 643 4390 643 4390 643
% 
\special{sh 1.000}%
\special{ia 4390 594 50 50  0.0000000  6.2831853}%
\special{pn 8}%
\special{ar 4390 594 50 50  0.0000000  6.2831853}%
% CIRCLE 2 0 0 0 Black Black
% 4 4290 993 4290 1043 4290 1043 4290 1043
% 
\special{sh 1.000}%
\special{ia 4290 994 50 50  0.0000000  6.2831853}%
\special{pn 8}%
\special{ar 4290 994 50 50  0.0000000  6.2831853}%
% LINE 2 2 3 0 Black Black
% 2 3990 593 4090 993
% 
\special{pn 8}%
\special{pa 3990 594}%
\special{pa 4090 994}%
\special{dt 0.045}%
% LINE 0 0 3 0 Black Black
% 2 4190 593 4290 993
% 
\special{pn 20}%
\special{pa 4190 594}%
\special{pa 4290 994}%
\special{fp}%
% LINE 0 0 3 0 Black Black
% 2 4390 593 4290 993
% 
\special{pn 20}%
\special{pa 4390 594}%
\special{pa 4290 994}%
\special{fp}%
% LINE 0 0 3 0 Black Black
% 2 4190 593 4090 993
% 
\special{pn 20}%
\special{pa 4190 594}%
\special{pa 4090 994}%
\special{fp}%
% SPLINE 2 2 3 0 Black Black
% 4 3600 593 3400 393 2800 503 2800 503
% 
\special{pn 8}%
\special{pn 8}%
\special{pa 3600 594}%
\special{pa 3595 587}%
\special{fp}%
\special{pa 3573 558}%
\special{pa 3568 551}%
\special{fp}%
\special{pa 3544 523}%
\special{pa 3538 517}%
\special{fp}%
\special{pa 3514 488}%
\special{pa 3509 483}%
\special{fp}%
\special{pa 3482 456}%
\special{pa 3476 450}%
\special{fp}%
\special{pa 3448 426}%
\special{pa 3442 421}%
\special{fp}%
\special{pa 3411 400}%
\special{pa 3404 396}%
\special{fp}%
\special{pa 3370 380}%
\special{pa 3363 377}%
\special{fp}%
\special{pa 3327 366}%
\special{pa 3319 365}%
\special{fp}%
\special{pa 3282 362}%
\special{pa 3274 361}%
\special{fp}%
\special{pa 3237 361}%
\special{pa 3229 362}%
\special{fp}%
\special{pa 3192 365}%
\special{pa 3184 367}%
\special{fp}%
\special{pa 3147 374}%
\special{pa 3139 375}%
\special{fp}%
\special{pa 3103 384}%
\special{pa 3095 386}%
\special{fp}%
\special{pa 3059 397}%
\special{pa 3052 400}%
\special{fp}%
\special{pa 3017 413}%
\special{pa 3009 416}%
\special{fp}%
\special{pa 2974 428}%
\special{pa 2966 431}%
\special{fp}%
\special{pa 2932 445}%
\special{pa 2924 448}%
\special{fp}%
\special{pa 2890 463}%
\special{pa 2883 467}%
\special{fp}%
\special{pa 2849 482}%
\special{pa 2841 486}%
\special{fp}%
\special{pa 2807 501}%
\special{pa 2800 504}%
\special{fp}%
% SPLINE 2 2 3 0 Black Black
% 6 4100 993 4100 593 4300 433 4700 433 4800 503 4800 503
% 
\special{pn 8}%
\special{pn 8}%
\special{pa 4100 994}%
\special{pa 4098 986}%
\special{fp}%
\special{pa 4091 949}%
\special{pa 4089 941}%
\special{fp}%
\special{pa 4082 904}%
\special{pa 4081 896}%
\special{fp}%
\special{pa 4075 858}%
\special{pa 4074 850}%
\special{fp}%
\special{pa 4069 813}%
\special{pa 4068 804}%
\special{fp}%
\special{pa 4066 767}%
\special{pa 4066 758}%
\special{fp}%
\special{pa 4066 721}%
\special{pa 4067 712}%
\special{fp}%
\special{pa 4072 675}%
\special{pa 4073 667}%
\special{fp}%
\special{pa 4083 631}%
\special{pa 4087 623}%
\special{fp}%
\special{pa 4103 589}%
\special{pa 4107 582}%
\special{fp}%
\special{pa 4130 552}%
\special{pa 4135 546}%
\special{fp}%
\special{pa 4162 518}%
\special{pa 4167 513}%
\special{fp}%
\special{pa 4197 489}%
\special{pa 4204 485}%
\special{fp}%
\special{pa 4237 465}%
\special{pa 4244 461}%
\special{fp}%
\special{pa 4277 444}%
\special{pa 4284 440}%
\special{fp}%
\special{pa 4320 427}%
\special{pa 4328 424}%
\special{fp}%
\special{pa 4363 411}%
\special{pa 4371 409}%
\special{fp}%
\special{pa 4408 400}%
\special{pa 4416 398}%
\special{fp}%
\special{pa 4453 393}%
\special{pa 4462 392}%
\special{fp}%
\special{pa 4499 390}%
\special{pa 4507 390}%
\special{fp}%
\special{pa 4545 390}%
\special{pa 4553 391}%
\special{fp}%
\special{pa 4591 396}%
\special{pa 4599 398}%
\special{fp}%
\special{pa 4636 406}%
\special{pa 4644 409}%
\special{fp}%
\special{pa 4679 423}%
\special{pa 4686 427}%
\special{fp}%
\special{pa 4719 445}%
\special{pa 4726 449}%
\special{fp}%
\special{pa 4757 471}%
\special{pa 4764 476}%
\special{fp}%
\special{pa 4794 499}%
\special{pa 4800 504}%
\special{fp}%
% STR 2 0 3 0 Black Black
% 4 3800 1093 3800 1193 5 0 0 0
% $A_{j}$
\put(38.0000,-11.9300){\makebox(0,0){$A_{j}$}}%
% STR 2 0 3 0 Black Black
% 4 3000 700 3000 800 5 0 0 0
% $Q'_{j}$
\put(30.0000,-8.0000){\makebox(0,0){$Q'_{j}$}}%
% STR 2 0 3 0 Black Black
% 4 4600 700 4600 800 5 0 0 0
% $Q''_{j}$
\put(46.0000,-8.0000){\makebox(0,0){$Q''_{j}$}}%
% STR 2 0 3 0 Black White
% 4 500 200 500 300 5 0 0 0
% $e_{j}$
\put(5.0000,-3.0000){\makebox(0,0){$e_{j}$}}%
% STR 2 0 3 0 Black White
% 4 3100 200 3100 300 5 0 0 0
% $e_{j}$
\put(31.0000,-3.0000){\makebox(0,0){$e_{j}$}}%
% STR 2 0 3 0 Black White
% 4 1900 200 1900 300 5 0 0 0
% $e_{j-1}$
\put(19.0000,-3.0000){\makebox(0,0){$e_{j-1}$}}%
% STR 2 0 3 0 Black White
% 4 4500 200 4500 300 5 0 0 0
% $e_{j-1}$
\put(45.0000,-3.0000){\makebox(0,0){$e_{j-1}$}}%
\end{picture}}%
\caption{Paths $Q'_{j}$ and $Q''_{j}$}
\label{f-s3-1}
\end{center}
\end{figure}
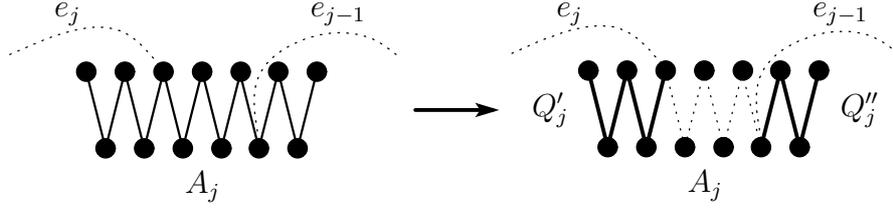

\begin{claim}%%%%%%%%%%%%%%%%%%%%%%%%%%%%%%%%%%%%%%%%%%%%%%%%%%%%%%%%%%%%%%%%%%%%%%%%%%%%%%%%%%%%%%%%%%%%%%%%%%%%%%%%%%%
\label{lemcl3.1.1.1}
Every path $A\in \A$ with $|V(A)|=2$ satisfies that $V(A)\cap T\subseteq T_{1}$.
\end{claim}
%%%%%%%%%%%%%%%%%%%%%%%%%%%%%%%%%%%%%%%%%%%%%%%%%%%%%%%%%%%%%%%%%%%%%%%%%%%%%%%%%%%%%%%%%%%%%%%%%%%%%%%%%%%%%%%%%%%%%%%%
\proof
Suppose that $\A$ contains a path $A$ such that $|V(A)|=2$ and $V(A)\cap T\not\subseteq T_{1}$ (i.e., $V(A)\cap T\subseteq T_{2}$).
Let $i$ be the minimum integer such that $\A_{i}$ contains a path $A_{i}$ such that $|V(A_{i})|=2$ and $V(A_{i})\cap T\subseteq T_{2}$.
Write $A_{i}=v^{(i)}_{1}v^{(i)}_{2}$, where $v^{(i)}_{1}\in S$ and $v^{(i)}_{2}\in T_{2}$, and set $l_{i}=2$.
By the minimality of $i$, every path $A$ belonging to $\bigcup _{1\leq j\leq i-1}\A_{j}$ with $|V(A)|=2$ satisfies $V(A)\cap T\subseteq T_{1}$.
By the definition of $\A_{j}$, there exist paths $A_{j}=v^{(j)}_{1}\cdots v^{(j)}_{l_{j}}\in \A_{j}~(0\leq j\leq i-1)$ such that $E_{H}(V(A_{j+1})\cap S,V(A_{j})\cap T)\not=\emptyset $ for every $j~(0\leq j\leq i-1)$.
For each $j~(0\leq j\leq i-1)$, we fix an edge $e_{j}\in E_{H}(V(A_{j+1})\cap S,V(A_{j})\cap T)$, and write $e_{j}=v^{(j+1)}_{s_{j+1}}v^{(j)}_{t_{j}}$.
By renumbering the vertices $v^{(j)}_{1},\ldots ,v^{(j)}_{l_{j}}$ of $A_{j}$ backward (i.e., by tracing the path $v^{(j)}_{1}\cdots v^{(j)}_{l_{j}}$ backward and numbering the vertices accordingly) if necessary, we may assume that $t_{j}<s_{j}$ for each $j~(1\leq j\leq i-1)$.
For each $j~(0\leq j\leq i-1)$, let $Q'_{j}$ be the path on $A_{j}$ from $v^{(j)}_{1}$ to $v^{(j)}_{t_{j}}$.
For each $j~(1\leq j\leq i)$, let $Q''_{j}$ be the path on $A_{j}$ from $v^{(j)}_{s_{j}}$ to $v^{(j)}_{l_{j}}$ (see Figure~\ref{f-s3-1}).
Note that if $A_{j}$ satisfies (II'), then $|V(Q'_{j})|$ is odd and $|V(Q''_{j})|$ is even.

Write $\{j\mid 1\leq j\leq i-1,|V(A_{j})|\geq 3\}=\{k_{1},k_{2},\ldots ,k_{m-1}\}$ with $1\leq k_{1}<k_{2}<\cdots <k_{m-1}\leq i-1$, and let $k_{0}=0$ and $k_{m}=i$ (it is possible that $m=1$).

\begin{figure}
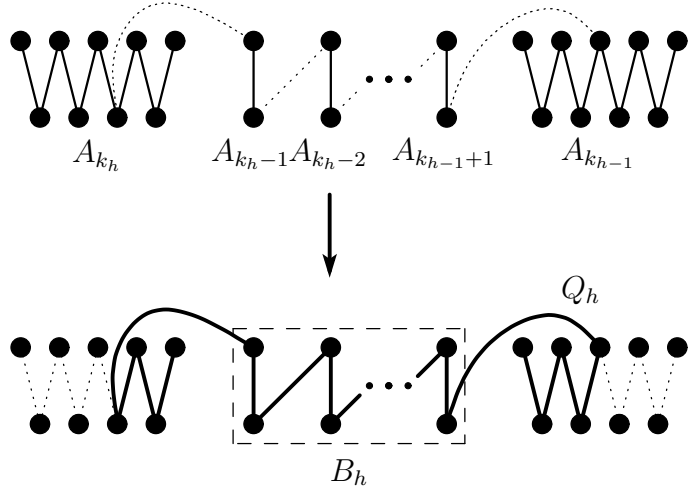

\begin{center}
%WinTpicVersion4.26
\unitlength 0.1in
%
\caption{Paths $B_{h}$ and $Q_{h}$}
\label{f-s3-2}
\end{center}
\end{figure}

Recall that every $A\in \bigcup _{1\leq j\leq i-1}\A_{j}$ with $|V(A)|=2$ satisfies $V(A)\cap T\subseteq T_{1}$.
Hence for each $h~(1\leq h\leq m)$, the graph $B_{h}=(\bigcup _{k_{h-1}+1\leq j\leq k_{h}-1}A_{j})+\{e_{j}\mid k_{h-1}+1\leq j\leq k_{h}-2\}$ is a path of $H$ with $V(B_{h})\cap T\subseteq T_{1}$ (here $B_{h}$ may be an empty graph).
Therefore for each $h~(1\leq h\leq m)$, the graph
$$
Q_{h}=(Q'_{k_{h-1}}\cup B_{h}\cup Q''_{k_{h}})+\{e_{k_{h-1}},e_{k_{h}-1}\}
$$
is a path of $H$ satisfying (II') (see Figure~\ref{f-s3-2}).
Note that when $h=m$, we here use the assumption that $V(A_{i})\cap T\subseteq T_{2}$.
Further, for $1\leq h\leq m-1$, since $|V(A_{k_{h}})|$ and $|V(Q'_{k_{h}})|$ are odd and $|V(Q''_{k_{h}})|$ is even, $A_{k_{h}}-(V(Q'_{k_{h}})\cup V(Q''_{k_{h}}))$ is a path of even order, and hence it has a perfect matching $M_{h}$.

Let
$$
F'=\left(F-\bigcup _{1\leq j\leq i}V(A_{j})\right)\cup \left(\bigcup _{1\leq h\leq m}Q_{h}\right)\cup \left(\bigcup _{1\leq h\leq m-1}M_{h}\right).
$$
Then $F'$ is a subgraph of $H$ such that $V(F')=V(F)\cup V(A_{0})~(=V(F)\cup \{v^{(0)}_{1}\})$ and each $A\in \C(F')$ is a path satisfying (I') or (II'), which contradicts the minimality of $|(S\cup T_{2})-V(F)|$, completing the proof of Claim~\ref{lemcl3.1.1.1}.
\qed

Let $Y_{0}=(\bigcup _{A\in \A}V(A))\cap T$.

\begin{claim}%%%%%%%%%%%%%%%%%%%%%%%%%%%%%%%%%%%%%%%%%%%%%%%%%%%%%%%%%%%%%%%%%%%%%%%%%%%%%%%%%%%%%%%%%%%%%%%%%%%%%%%%%%%
\label{lemcl3.1.1.2}
We have $N_{H}(Y_{0})=(\bigcup _{A\in \A}V(A))\cap S$.
\end{claim}
%%%%%%%%%%%%%%%%%%%%%%%%%%%%%%%%%%%%%%%%%%%%%%%%%%%%%%%%%%%%%%%%%%%%%%%%%%%%%%%%%%%%%%%%%%%%%%%%%%%%%%%%%%%%%%%%%%%%%%%%
\proof
Suppose that $N_{H}(Y_{0})\not=(\bigcup _{A\in \A}V(A))\cap S$.
Then there exists an integer $i$ and there exists a vertex $v\in S-(\bigcup _{A\in \A}V(A))$ such that $N_{H}(v)\cap (\bigcup _{A\in \A_{i}}V(A))\not=\emptyset $.
Let $A_{i+1}$ be the path of $H$ consisting of $v$.
By the definition of $\A_{j}$, there exist paths $A_{j}\in \A_{j}~(0\leq j\leq i)$ such that $E_{H}(V(A_{j+1})\cap S,V(A_{j})\cap T)\not=\emptyset $ for every $j~(0\leq j\leq i)$.
For each $j~(0\leq j\leq i)$, we fix an edge $u_{j}v_{j+1}\in E_{H}(V(A_{j+1})\cap S,V(A_{j})\cap T)$ with $u_{j}\in V(A_{j})\cap T$ and $v_{j+1}\in V(A_{j+1})\cap S$.

Let $k~(0\leq k\leq i)$ be the maximum integer such that $|V(A_{k})|$ is odd (the fact that $|V(A_{0})|=1$ assures us the existence of $k$).
Then for each $j~(k+1\leq j\leq i)$, we have $|V(A_{j})|=2$ (i.e., $A_{j}=u_{j}v_{j}$).
Furthermore, since $A_{k}$ is a path with $|V(A_{k})\cap T|=|V(A_{k})\cap S|+1$ and $u_{k}\in T$, $A_{k}-u_{k}$ has a perfect matching $M$.
Hence $M^{*}=\{u_{j}v_{j+1}\mid k\leq j\leq i\}\cup M$ is a perfect matching of the subgraph of $H$ induced by $\bigcup _{k\leq j\leq i+1}V(A_{j})$.
Therefore $F'=(F-\bigcup _{k\leq j\leq i}V(A_{j}))\cup M^{*}$ is a subgraph of $H$ such that $V(F')\supseteq V(F)\cup \{v\}$ and each $A\in \C(F')$ is a path satisfying (I') or (II'), which contradicts the minimality of $|(S\cup T_{2})-V(F)|$.
\qed

We continue with the proof of the lemma.
By the definition of $\A$, we have
\begin{align}
\label{eq-sec3-1}
Y_{0}\cap T_{1}=\left(\bigcup _{A\in \A-\A_{0}}V(A)\right)\cap T_{1}
\end{align}
and
\begin{align}
\label{eq-sec3-2}
Y_{0}\cap T_{2}=\left(\left(\bigcup _{A\in \A-\A_{0}}V(A)\right)\cap T_{2}\right)\cup (T_{2}-V(F)).
\end{align}
If $A\in \A$ satisfies (I'), then $|V(A)\cap S|=1=|V(A)\cap T_{1}|$ and $V(A)\cap T_{2}=\emptyset $ by Claim~\ref{lemcl3.1.1.1}.
Thus
\begin{align}
\label{eq-sec3-3}
|V(A)\cap S|=|V(A)\cap T_{1}|+\frac{1}{2}|V(A)\cap T_{2}|\mbox{~~~for each }A\in \A \mbox{ satisfying (I')}.
\end{align}
If $A\in \A$ satisfies (II'), then $|V(A)\cap T_{1}|=|V(A)\cap S|-1$ and $|V(A)\cap T_{2}|=2$ by (II').
Thus
\begin{align}
\label{eq-sec3-4}
|V(A)\cap S|=|V(A)\cap T_{1}|+\frac{1}{2}|V(A)\cap T_{2}|\mbox{~~~for each }A\in \A \mbox{ satisfying (II')}.
\end{align}
Recall that $T_{2}-V(F)\not=\emptyset $.
Hence by Claim~\ref{lemcl3.1.1.2} and (\ref{eq-sec3-1})--(\ref{eq-sec3-4}),
\begin{align*}
|N_{H}(Y_{0})| &= \sum _{A\in \A}|V(A)\cap S|\\
&= \sum _{A\in \A-\A_{0}}|V(A)\cap S|\\
&= \sum _{A\in \A-\A_{0}}\left(|V(A)\cap T_{1}|+\frac{1}{2}|V(A)\cap T_{2}|\right)\\
&= \sum _{A\in \A-\A_{0}}|V(A)\cap T_{1}|+\frac{1}{2}\sum _{A\in \A-\A_{0}}|V(A)\cap T_{2}|\\
&= |Y_{0}\cap T_{1}|+\frac{1}{2}(|Y_{0}\cap T_{2}|-|T_{2}-V(F)|)\\
&< |Y_{0}\cap T_{1}|+\frac{1}{2}|Y_{0}\cap T_{2}|,
\end{align*}
which contradicts the assumption of the lemma.

This completes the proof of Lemma~\ref{lem3.1.1}.
\qed

\medbreak\noindent\textit{Proof of Proposition~\ref{prop3.1}.}\quad
Applying Lemma~\ref{lem3.1.1} to $(G-L)[S\cup T_{1}\cup T_{2}]$, we see that $G-L$ has a subgraph $F'$ with $V(F')\cap T=T_{1}\cup T_{2}$ such that each $A\in \C(F')$ is a path with $V(A)\cap T\subseteq T_{1}\cup T_{2}$ satisfying (I) or (II).
In particular, $G$ has a subgraph $F$ with $V(F)\supseteq T_{1}\cup T_{2}$ such that each $A\in \C(F)$ is a path satisfying (I) or (II).
Choose $F$ so that $|S-V(F)|$ is as small as possible.

It suffices to show that $S-V(F)=\emptyset $.
By way of contradiction, suppose that $S-V(F)\not=\emptyset $.
Now we define the set $\A$ of paths of $G$ as follows:
Let $\A_{0}$ be the set of paths of $G$ consisting of one vertex in $S-V(F)$.
Let $\D$ be the set of paths of $G$ consisting of one vertex in $T-V(F)$.
For each $i\geq 1$, let $\A_{i}$ be the set of those members $A$ of $\C(F)\cup \D$ such that $A\not\in \bigcup _{0\leq j\leq i-1}\A_{j}$ and $E_{G}(V(A)\cap T,\bigcup _{A'\in \A_{i-1}}(V(A')\cap S))\not=\emptyset $.
Set $\A=\bigcup _{i\geq 0}\A_{i}$.

Suppose that $\A-\A_{0}$ contains a path of odd order.
Let $i$ be the minimum integer such that $\A_{i}$ contains a path $A_{i}$ of odd order.
By the definition of $\A_{j}$, there exist paths $A_{j}\in \A_{j}~(0\leq j\leq i-1)$ such that $E_{G}(V(A_{j+1})\cap T,V(A_{j})\cap S)\not=\emptyset $ for every $0\leq j\leq i-1$.
Write $V(A_{0})=\{v_{0}\}$.
By the minimality of $i$, for each $j~(1\leq j\leq i-1)$, we have $|V(A_{j})|=2$.
For each $j~(1\leq j\leq i-1)$, write $A_{j}=u_{j}v_{j}$, where $V(A_{j})\cap T=\{u_{j}\}$ and $V(A_{j})\cap S=\{v_{j}\}$.
Let $u_{i}\in N_{G}(v_{i-1})\cap V(A_{i})$.
%For each $j~(0\leq j\leq i-1)$, we fix an edge $v_{j}u_{j+1}\in E_{G}(V(A_{j+1})\cap T,V(A_{j})\cap S)$ with $v_{j}\in V(A_{j})\cap S$ and $u_{j+1}\in V(A_{j+1})\cap T$.
Since $A_{i}$ is a path with $|V(A_{i})\cap T|=|V(A_{i})\cap S|+1$ and $u_{i}\in T$, $A_{i}-u_{i}$ has a perfect matching $M$.
Hence $M^{*}=\{v_{j}u_{j+1}\mid 0\leq j\leq i-1\}\cup M$ is a perfect matching of the subgraph of $G$ induced by $\bigcup _{0\leq j\leq i}V(A_{j})$.
Therefore $F'=(F-\bigcup _{1\leq j\leq i}V(A_{j}))\cup M^{*}$ is a subgraph of $G$ such that $V(F')\supseteq V(F)\cup \{v_{0}\}$ and each $A\in \C(F')$ is a path satisfying (I) or (II), which contradicts the minimality of $|S-V(F)|$.
Thus every element of $\A-\A_{0}$ is a path of order $2$.
In particular, $\A\cap \D=\emptyset $.

Let $X_{0}=(\bigcup _{A\in \A}V(A))\cap S$.
Since $\A\cap \D=\emptyset $, $N_{G}(X_{0})=(\bigcup _{A\in \A-\A_{0}}V(A))\cap T$.
Since every element of $\A-\A_{0}$ is a path of order $2$, $|(\bigcup _{A\in \A-\A_{0}}V(A))\cap T|=|(\bigcup _{A\in \A-\A_{0}}V(A))\cap S|$.
Consequently
\begin{align*}
|N_{G}(X_{0})| &= \sum _{A\in \A-\A_{0}}|V(A)\cap T|\\
&= \sum _{A\in \A-\A_{0}}|V(A)\cap S|\\
&= \sum _{A\in \A}|V(A)\cap S|-|S-V(F)|\\
&< \sum _{A\in \A}|V(A)\cap S|\\
&= |X_{0}|,
\end{align*}
which contradicts the assumption of the proposition.
\qed

%%%%%%%%%%%%%%%%%%%%%%%%%%%%%%%%%%%%%%%%%%%%%%%%%%%%%%%%%%%%%%%%%%%%%%%%%%%%%%%%%%%%%%%%%%%%%%%%%%%%%%%%%%%%%%%%%%%%%%%%
%%%%%%%%%%%%%%%%%%%%%%%%%%%%%%%%%%%%%%%%%%%%%%%%%%%%%%%%%%%%%%%%%%%%%%%%%%%%%%%%%%%%%%%%%%%%%%%%%%%%%%%%%%%%%%%%%%%%%%%%
%%%%%%%%%%%%%%%%%%%%%%%%%%%%%%%%%%%%%%%%%%%%%%%%%%%%%%%%%%%%%%%%%%%%%%%%%%%%%%%%%%%%%%%%%%%%%%%%%%%%%%%%%%%%%%%%%%%%%%%%
\section{Hypomatchable graphs having no $\{P_{2},P_{2k+1}\}$-factor}\label{sec4}
%%%%%%%%%%%%%%%%%%%%%%%%%%%%%%%%%%%%%%%%%%%%%%%%%%%%%%%%%%%%%%%%%%%%%%%%%%%%%%%%%%%%%%%%%%%%%%%%%%%%%%%%%%%%%%%%%%%%%%%%
%%%%%%%%%%%%%%%%%%%%%%%%%%%%%%%%%%%%%%%%%%%%%%%%%%%%%%%%%%%%%%%%%%%%%%%%%%%%%%%%%%%%%%%%%%%%%%%%%%%%%%%%%%%%%%%%%%%%%%%%
%%%%%%%%%%%%%%%%%%%%%%%%%%%%%%%%%%%%%%%%%%%%%%%%%%%%%%%%%%%%%%%%%%%%%%%%%%%%%%%%%%%%%%%%%%%%%%%%%%%%%%%%%%%%%%%%%%%%%%%%

A graph $G$ is {\it hypomatchable} if $G-x$ has a perfect matching for every $x\in V(G)$.
In this section, we characterize hypomatchable graphs having no $\{P_{2},P_{2k+1}\}$-factor for $k\in \{3,4\}$.

%%%%%%%%%%%%%%%%%%%%%%%%%%%%%%%%%%%%%%%%%%%%%%%%%%%%%%%%%%%%%%%%%%%%%%%%%%%%%%%%%%%%%%%%%%%%%%%%%%%%%%%%%%%%%%%%%%%%%%%%
%%%%%%%%%%%%%%%%%%%%%%%%%%%%%%%%%%%%%%%%%%%%%%%%%%%%%%%%%%%%%%%%%%%%%%%%%%%%%%%%%%%%%%%%%%%%%%%%%%%%%%%%%%%%%%%%%%%%%%%%
\subsection{Fundamental properties of hypomatchable graphs}\label{sec4.1}
%%%%%%%%%%%%%%%%%%%%%%%%%%%%%%%%%%%%%%%%%%%%%%%%%%%%%%%%%%%%%%%%%%%%%%%%%%%%%%%%%%%%%%%%%%%%%%%%%%%%%%%%%%%%%%%%%%%%%%%%
%%%%%%%%%%%%%%%%%%%%%%%%%%%%%%%%%%%%%%%%%%%%%%%%%%%%%%%%%%%%%%%%%%%%%%%%%%%%%%%%%%%%%%%%%%%%%%%%%%%%%%%%%%%%%%%%%%%%%%%%

We start with a structure theorem for hypomatchable graphs.
Let $G$ be a graph.
A sequence $(H_{1},\ldots ,H_{m})$ of edge-disjoint subgraphs of $G$ is an {\it ear decomposition} if
\begin{enumerate}[{\bf (E1)}]
\item
$V(G)=\bigcup _{1\leq i\leq m}V(H_{i})$;
\item
for each $1\leq i\leq m$, $|E(H_{i})|$ is odd and $|E(H_{i})|\geq 3$;
\item
$H_{1}$ is a cycle; and
\item
for each $2\leq i\leq m$, either
\begin{enumerate}
\item[{\bf (E4-1)}]
$H_{i}$ is a path and only the endvertices of $H_{i}$ belong to $\bigcup _{1\leq j\leq i-1}V(H_{j})$, or
\item[{\bf (E4-2)}]
$H_{i}$ is a cycle with $|V(H_{i})\cap (\bigcup _{1\leq j\leq i-1}V(H_{j}))|=1$.
\end{enumerate}
\end{enumerate}

Lov\'{a}sz~\cite{L} proved the following theorem.

\begin{Thm}[Lov\'{a}sz~\cite{L}]%%%%%%%%%%%%%%%%%%%%%%%%%%%%%%%%%%%%%%%%%%%%%%%%%%%%%%%%%%%%%%%%%%%%%%%%%%%%%%%%%%%%%%%%
\label{Thm4-1.A}
Let $G$ be a graph with $|V(G)|\geq 3$.
\begin{enumerate}[{\upshape(i)}]
\item
If $G$ has an ear decomposition, then $G$ is hypomatchable.
\item
If $G$ is hypomatchable, then for each $e\in E(G)$, $G$ has an ear decomposition $(H_{1},\ldots ,H_{m})$ such that $e\in E(H_{1})$.
\end{enumerate}
\end{Thm}
%%%%%%%%%%%%%%%%%%%%%%%%%%%%%%%%%%%%%%%%%%%%%%%%%%%%%%%%%%%%%%%%%%%%%%%%%%%%%%%%%%%%%%%%%%%%%%%%%%%%%%%%%%%%%%%%%%%%%%%%

In the remainder of this subsection, we let $G$ be a hypomatchable graph, and let $\H=(H_{1},\ldots ,H_{m})$ be an ear decomposition of $G$.
We start with lemmas which hold for an ear decomposition of a hypomatchable graph in general.

\begin{lem}%%%%%%%%%%%%%%%%%%%%%%%%%%%%%%%%%%%%%%%%%%%%%%%%%%%%%%%%%%%%%%%%%%%%%%%%%%%%%%%%%%%%%%%%%%%%%%%%%%%%%%%%%%%%%
\label{lem4-1.0}
For each $i~(2\leq i\leq m)$, there exists an ear decomposition $(H'_{1},\ldots ,H'_{m'})$ of $G$ such that $H_{i}\subseteq H'_{1}$.
\end{lem}
%%%%%%%%%%%%%%%%%%%%%%%%%%%%%%%%%%%%%%%%%%%%%%%%%%%%%%%%%%%%%%%%%%%%%%%%%%%%%%%%%%%%%%%%%%%%%%%%%%%%%%%%%%%%%%%%%%%%%%%%
\proof
Set $H=H_{1}\cup \cdots \cup H_{i}$.
Then $(H_{1},\ldots ,H_{i})$ is an ear decomposition of $H$, and hence $H$ is hypomatchable by Theorem~\ref{Thm4-1.A}(i).
Take $e\in E(H_{i})$.
By Theorem~\ref{Thm4-1.A}(ii), $H$ has an ear decomposition $(H'_{1},\ldots ,H'_{n})$ such that $e\in E(H'_{1})$.
Since $H_{i}$ satisfies (E4), we have $d_{H}(v)=2$ for all $v\in V(H_{i})-(\bigcup _{1\leq j\leq i-1}V(H_{j}))$.
Since $H'_{1}$ satisfies (E3), this implies $H_{i}\subseteq H'_{1}$.
Since $\bigcup _{1\leq j\leq n}V(H'_{j})=\bigcup _{1\leq j\leq i}V(H_{j})$, it follows that $(H'_{1},\ldots ,H'_{n},H_{i+1},\ldots ,H_{m})$ is an ear decomposition of $G$ with the desired property.
\qed

\begin{lem}%%%%%%%%%%%%%%%%%%%%%%%%%%%%%%%%%%%%%%%%%%%%%%%%%%%%%%%%%%%%%%%%%%%%%%%%%%%%%%%%%%%%%%%%%%%%%%%%%%%%%%%%%%%%%
\label{lem4-1.00}
Suppose that each $H_{i}~(1\leq i\leq m)$ is a cycle, and let $i_{1},\ldots ,i_{m}$ be a permutation of $1,\ldots ,m$ such that $V(H_{i_{l}})\cap (\bigcup _{1\leq j\leq l-1}V(H_{i_{j}}))\not=\emptyset $ for each $l~(2\leq l\leq m)$.
Then $(H_{i_{1}},\ldots ,H_{i_{m}})$  is an ear decomposition of $G$.
\end{lem}
%%%%%%%%%%%%%%%%%%%%%%%%%%%%%%%%%%%%%%%%%%%%%%%%%%%%%%%%%%%%%%%%%%%%%%%%%%%%%%%%%%%%%%%%%%%%%%%%%%%%%%%%%%%%%%%%%%%%%%%%
\proof
Since each $H_{i}$ is a cycle, it follows from the definition of an ear decomposition that $H_{i}$ is a block of $G$ for each $i$.
Thus for each $l~(2\leq l\leq m)$, the assumption that $V(H_{i_{l}})\cap (\bigcup _{1\leq j\leq l-1}V(H_{i_{j}}))\not=\emptyset $ implies that $|V(H_{i_{l}})\cap (\bigcup _{1\leq j\leq l-1}V(H_{i_{j}}))|=1$.
Hence by the definition of an ear decomposition, $(H_{i_{1}},\ldots ,H_{i_{m}})$ is also an ear decomposition.
\qed

Our next result is concerned with a hypomatchable graph with no $\{P_{2},P_{2k+1}\}$-factor.
In order to state the result, we need some more definitions.
For each $i~(1\leq i\leq m)$, let $P_{\H}(i)=H_{i}-\bigcup _{1\leq j\leq i-1}V(H_{j})$.
Note that $V(P_{\H}(i))\cap V(H_{j})=\emptyset $ for any $i,j$ with $i>j$, and $\bigcup _{1\leq j\leq i}V(H_{j})=\bigcup _{1\leq j\leq i}V(P_{\H}(j))$ for each $i$.
We have $P_{\H}(1)=H_{1}$ and, by (E2) and (E4), $P_{\H}(i)$ is a path of even order for $2\leq i\leq m$.
For an odd integer $s\geq 5$, a set $I\subseteq \{1,2,\ldots ,m\}$ of indices with $1\in I$ is {\it $s$-large} with respect to $\H$ if $\sum _{i\in I}|V(P_{\H}(i))|\geq s$ and the subgraph of $G$ induced by $\bigcup _{i\in I}V(P_{\H}(i))$ has a spanning path.

\begin{lem}%%%%%%%%%%%%%%%%%%%%%%%%%%%%%%%%%%%%%%%%%%%%%%%%%%%%%%%%%%%%%%%%%%%%%%%%%%%%%%%%%%%%%%%%%%%%%%%%%%%%%%%%%%%%%
\label{lem4-1.1}
Let $k\geq 3$, and suppose that $G$ has no $\{P_{2},P_{2k+1}\}$-factor.
Then there is no $(2k+1)$-large set with respect to $\H$.
\end{lem}
%%%%%%%%%%%%%%%%%%%%%%%%%%%%%%%%%%%%%%%%%%%%%%%%%%%%%%%%%%%%%%%%%%%%%%%%%%%%%%%%%%%%%%%%%%%%%%%%%%%%%%%%%%%%%%%%%%%%%%%%
\proof
Suppose that there exists a $(2k+1)$-large set $I$ with respect to $\H$.
Then by Fact~\ref{fact1}, the subgraph of $G$ induced by $\bigcup _{i\in I}V(P_{\H}(i))$ has a $\{P_{2},P_{2k+1}\}$-factor $F$.
On the other hand, for each $i$ with $2\leq i\leq m$ and $i\not\in I$, from the fact that $P_{\H}(i)$ is a path of even order, we see that $P_{\H}(i)$ has a perfect matching $M_{i}$.
Since $\{V(P_{\H}(i))\mid i\not\in I\}$ is a partition of $V(G)-(\bigcup _{i\in I}V(P_{\H}(i)))$, $F\cup (\bigcup _{i\not\in I}M_{i})$ is a $\{P_{2},P_{2k+1}\}$-factor of $G$, which is a contradiction.
\qed

Throughout the rest of this subsection, we assume that we have chosen $\H=(H_{1},\ldots ,H_{m})$ so that
\begin{enumerate}[{\bf (H1)}]
\item
$|E(H_{1})|$ is as large as possible.
\end{enumerate}

\begin{lem}%%%%%%%%%%%%%%%%%%%%%%%%%%%%%%%%%%%%%%%%%%%%%%%%%%%%%%%%%%%%%%%%%%%%%%%%%%%%%%%%%%%%%%%%%%%%%%%%%%%%%%%%%%%%%
\label{lem4-1.20}
Suppose that $\H=(H_{1},\ldots ,H_{m})$ is chosen so that (H1) holds.
Let $2\leq i\leq m$, and let $v$, $v'$ be the endvertices of $P_{\H}(i)$.
Then no two vertices $w$, $w'$ with $w\in N_{G}(v)\cap V(H_{1})$ and $w'\in N_{G}(v')\cap V(H_{1})$ are consecutive on $H_{1}$.
\end{lem}
%%%%%%%%%%%%%%%%%%%%%%%%%%%%%%%%%%%%%%%%%%%%%%%%%%%%%%%%%%%%%%%%%%%%%%%%%%%%%%%%%%%%%%%%%%%%%%%%%%%%%%%%%%%%%%%%%%%%%%%%
\proof
Suppose that there exist $w\in N_{G}(v)\cap V(H_{1})$ and $w'\in N_{G}(v')\cap V(H_{1})$ such that $w$ and $w'$ are consecutive on $H_{1}$.
Then $G[V(H_{1})\cup V(P_{\H}(i))]$ contains a spanning cycle $C$.
Since $|E(C)|=|V(C)|=|V(H_{1})|+|V(P_{\H}(i))|$, $|E(C)|$ is odd and $|E(C)|>|E(H_{1})|$.
Since $V(H_{j})\cap V(P_{\H}(i))=\emptyset $ for every $j$ with $2\leq j\leq i-1$, $(C,H_{2},\ldots ,H_{i-1})$ is an ear decomposition of $G[V(H_{1})\cup \cdots \cup V(H_{i})]$, and hence $(C,H_{2},\ldots ,H_{i-1},H_{i+1},\ldots ,H_{m})$ is an ear decomposition of $G$, which contradicts (H1).
\qed

\begin{lem}%%%%%%%%%%%%%%%%%%%%%%%%%%%%%%%%%%%%%%%%%%%%%%%%%%%%%%%%%%%%%%%%%%%%%%%%%%%%%%%%%%%%%%%%%%%%%%%%%%%%%%%%%%%%%
\label{lem4-1.21}
Suppose that (H1) holds, and suppose further that $|E(H_{1})|=3$. 
Then each $H_{i}$ is a cycle of order $3$, and $G=H_{1}\cup \cdots \cup H_{m}$.
\end{lem}
%%%%%%%%%%%%%%%%%%%%%%%%%%%%%%%%%%%%%%%%%%%%%%%%%%%%%%%%%%%%%%%%%%%%%%%%%%%%%%%%%%%%%%%%%%%%%%%%%%%%%%%%%%%%%%%%%%%%%%%%
\proof
Let $2\leq i\leq m$.
By Lemma~\ref{lem4-1.0}, there is an ear decomposition $(H'_{1},\ldots ,H'_{m'})$ such that $H_{i}\subseteq H'_{1}$.
If $|E(H_{i})|>3$ or $H_{i}$ is a path, then we get $|E(H'_{1})|>3$, which contradicts (H1).
Thus each $H_{i}$ is a cycle of order $3$.

Now suppose that there exists $e=ab\in E(G)$ such that $e\not\in E(H_{1}\cup \cdots \cup H_{m})$.
Since $(H_{1}\cup \cdots \cup H_{m})+e$ is hypomatchable by Theorem~\ref{Thm4-1.A}(i), it follows from Theorem~\ref{Thm4-1.A}(ii) that there is an ear decomposition $(H'_{1},\ldots ,H'_{m'})$ of $(H_{1}\cup \cdots \cup H_{m})+e$ such that $e\in E(H'_{1})$.
By (H1), $|E(H'_{1})|=3$.
Write $H'_{1}=abva$.
Let $i$, $j$ be the indices such that $av\in E(H_{i})$ and $bv\in E(H_{j})$.
Then $i\not=j$, $v\in V(H_{i})\cap V(H_{j})$, and $(H_{i}\cup H_{j})+e$ has a spanning cycle $C$.
By Lemma~\ref{lem4-1.00}, $G$ has an ear decomposition $(H''_{1},\ldots ,H''_{m})$ with $H''_{1}=H_{i}$ and $H''_{2}=H_{j}$.
This implies that $(C,H''_{3},\ldots ,H''_{m})$ is an ear decomposition of $G$, which contradicts (H1).
Thus $G=H_{1}\cup \cdots \cup H_{m}$.
\qed

%%%%%%%%%%%%%%%%%%%%%%%%%%%%%%%%%%%%%%%%%%%%%%%%%%%%%%%%%%%%%%%%%%%%%%%%%%%%%%%%%%%%%%%%%%%%%%%%%%%%%%%%%%%%%%%%%%%%%%%%
%%%%%%%%%%%%%%%%%%%%%%%%%%%%%%%%%%%%%%%%%%%%%%%%%%%%%%%%%%%%%%%%%%%%%%%%%%%%%%%%%%%%%%%%%%%%%%%%%%%%%%%%%%%%%%%%%%%%%%%%
\subsection{Constructions of hypomatchable graphs}\label{sec4.2}
%%%%%%%%%%%%%%%%%%%%%%%%%%%%%%%%%%%%%%%%%%%%%%%%%%%%%%%%%%%%%%%%%%%%%%%%%%%%%%%%%%%%%%%%%%%%%%%%%%%%%%%%%%%%%%%%%%%%%%%%
%%%%%%%%%%%%%%%%%%%%%%%%%%%%%%%%%%%%%%%%%%%%%%%%%%%%%%%%%%%%%%%%%%%%%%%%%%%%%%%%%%%%%%%%%%%%%%%%%%%%%%%%%%%%%%%%%%%%%%%%

In this subsection, we constructs five families $\G_{0},\G_{1},\G_{2},\G_{3},\G_{4}$ of hypomatchable graphs (see Figure~\ref{f-s4.2-1}).

\begin{figure}
\begin{center}
%WinTpicVersion4.28b
{\unitlength 0.1in
}%
\caption{Graphs in $\G^{*}_{i}$ or $\G_{i}$}
\label{f-s4.2-1}
\end{center}
\end{figure}

\begin{enumerate}[$\bullet $]
\item
Let $\G^{*}_{0}=\{K_{1}+sK_{2}\mid s\geq 2\}$ and $\G_{0}=\{K_{1}+sK_{2}\mid s\geq 3\}$.
Note that for each $H\in \G^{*}_{0}$, $H$ is hypomatchable and has no $\{P_{2},P_{7}\}$-factor.
\end{enumerate}

Let $s_{1},s_{2},s_{3}$ be nonnegative integers.
Let $Q=u_{1}u_{2}u_{3}$ be a path of order $3$ and, for $i\in \{1,2,3\}$ and $1\leq j\leq s_{i}$, let $L_{i,j}$ be a path of order $2$.
For each $1\leq j\leq s_{2}$, write $L_{2,j}=v_{1,j}v_{3,j}$.

\begin{enumerate}[$\bullet $]
\item
Let $A_{1}(s_{1},s_{2},s_{3})$ be the graph obtained from $Q\cup (\bigcup _{i\in \{1,2,3\}}(\bigcup _{1\leq j\leq s_{i}}L_{i,j}))$ by adding the edge $u_{1}u_{3}$ and joining $u_{i}$ to all vertices in $\bigcup _{1\leq j\leq s_{i}}V(L_{i,j})$ for each $i\in \{1,2,3\}$.
Note that $A_{1}(s_{1},0,0)\simeq K_{1}+(s_{1}+1)K_{2}$.
Let $\G^{*}_{1}=\{A_{1}(s_{1},s_{2},s_{3})\mid s_{1}+s_{2}+s_{3}\geq 1\}$ and $\G_{1}=\{A_{1}(s_{1},s_{2},s_{3})\mid s_{1}+s_{2}+s_{3}\geq 3\}$.

We divide the set $\G_{1}$ into three sets.
Let $\G^{(1)}_{1}=\{A_{1}(s_{1},s_{2},s_{3})\in \G_{1}\mid \min \{s_{1},s_{2},s_{3}\}\leq 1\}$, $\G^{(2)}_{1}=\{A_{1}(s_{1},s_{2},s_{3})\in \G_{1}\mid \min \{s_{1},s_{2},s_{3}\}=2\}$ and $\G^{(3)}_{1}=\{A_{1}(s_{1},s_{2},s_{3})\in \G_{1}\mid \min \{s_{1},s_{2},s_{3}\}\geq 3\}$.
\item
Let $A'_{2}(s_{1},s_{2},s_{3})$ be the graph obtained from $Q\cup (\bigcup _{i\in \{1,2,3\}}(\bigcup _{1\leq j\leq s_{i}}L_{i,j}))$ by joining $u_{i}$ to all vertices in $(\bigcup _{1\leq j\leq s_{i}}V(L_{i,j}))\cup \{v_{i,j}\mid 1\leq j\leq s_{2}\}$ for each $i\in \{1,3\}$.
Let $A''_{2}(s_{1},s_{2},s_{3})$ be the graph obtained from $Q\cup (\bigcup _{i\in \{1,2,3\}}(\bigcup _{1\leq j\leq s_{i}}L_{i,j}))$ by adding the edge $u_{1}u_{3}$ and joining $u_{i}$ to all vertices in $(\bigcup _{1\leq j\leq s_{i}}V(L_{i,j}))\cup (\bigcup _{1\leq j\leq s_{2}}V(L_{2,j}))$ for each $i\in \{1,3\}$.
Let $\G_{2}=\{H\mid A'_{2}(s_{1},s_{2},s_{3})\subseteq H\subseteq A''_{2}(s_{1},s_{2},s_{3})$ with $s_{2}\geq 1$, and $s_{1}+s_{2}+s_{3}\geq 3$, and either $s_{1}\geq 1$ and $s_{3}\geq 1$ or $s_{2}\geq 2\}$.
\item
Assume $s_{2}=1$ and $s_{3}=0$.
Let $A'_{3}(s_{1})=A'_{2}(s_{1},1,0)$.
Let $A''_{3}(s_{1})$ be the graph obtained from $A'_{3}(s_{1})$ by joining all possible pairs of vertices in $V(Q)\cup L_{2,1}$.
Let $\G_{3}=\{H\mid A'_{3}(s_{1})\subseteq H\subseteq A''_{3}(s_{1})$ with $s_{1}\geq 2\}$.
\item
Assume that $s_{2}=2$ and $s_{3}=0$.
Let $A'_{4}(s_{1})$ be the graph obtained from $A'_{2}(s_{1},2,0)$ by adding the edge $v_{3,1}v_{3,2}$.
Let $A''_{4}(s_{1})$ be the graph obtained from $A'_{4}(s_{1})$ by adding the edges $u_{1}u_{3},u_{1}v_{3,1},u_{1}v_{3,2}$.
Let $\G_{4}=\{H\mid A'_{4}(s_{1})\subseteq H\subseteq A''_{4}(s_{1})$ with $s_{1}\geq 1\}$.

We can verify that for each $H\in \G^{*}_{1}\cup \G_{2}\cup \G_{3}\cup \G_{4}$, $H$ is hypomatchable and has no $\{P_{2},P_{9}\}$-factor.
\end{enumerate}

Now we define crush sets of graphs belonging to $\bigcup _{0\leq i\leq 4}\G_{i}$.
For $H\in \G_{0}$, a set $X\subseteq V(H)$ is a {\it crush set} of $H$ if $x\in X$ and $|X\cap V(C)|=1$ for each $C\in \C(H-x)$, where $x$ is the unique cutvertex of $H$.
Let $H\in \G_{1}$, and write $H=A_{1}(s_{1},s_{2},s_{3})$.
We may assume that $\min\{s_{1},s_{2},s_{3}\}=s_{3}$.
If $H\in \G^{(1)}_{1}\cup \G^{(2)}_{1}$, a {\it crush set} of $H$ is a set $X\subseteq V(G)$ such that $X\cap V(Q)=\{u_{1},u_{2}\}$, $X\cap (\bigcup _{1\leq j\leq s_{3}}V(L_{3,j}))=\emptyset $ and $|X\cap V(L_{i,j})|=1$ for each $i\in \{1,2\}$ and each $1\leq j\leq s_{i}$ (note that if $s_{3}=0$ and $s_{1}$ or $s_{2}$ is zero, then this definition is consistent with the definition of a crush set for a graph in $\G_{0}$).
If $H\in \G^{(3)}_{1}$, a {\it crush set} of $H$ is a set $X\subseteq V(G)$ such that $V(Q)\subseteq X$ and $|X\cap V(L_{i,j})|=1$ for each $i\in \{1,2,3\}$ and each $1\leq j\leq s_{i}$.
For $H\in \G_{2}$, a set $X\subseteq V(H)$ is a {\it crush set} of $H$ if $X\cap V(Q)=\{u_{1},u_{3}\}$ and $|X\cap V(L_{i,j})|=1$ for each $i\in \{1,2,3\}$ and each $1\leq j\leq s_{i}$.
For $H\in \G_{3}$, a set $X\subseteq V(H)$ is a {\it crush set} of $H$ if $X\cap V(Q)=\{u_{1},u_{3}\}$, $X\cap V(L_{2,1})=\emptyset $ and $|X\cap V(L_{1,j})|=1$ for each $1\leq j\leq s_{1}$.
For $H\in \G_{4}$, a set $X\subseteq V(H)$ is a {\it crush set} of $H$ if $X\cap V(Q)=\{u_{1},u_{3}\}$, $X\cap V(L_{2,1})=\{v_{3,1}\}$, $X\cap V(L_{2,2})=\{v_{3,2}\}$ and $|X\cap V(L_{1,j})|=1$ for each $1\leq j\leq s_{1}$.

By inspection, we get the following lemma, which will be used in Section~\ref{sec5}.

\begin{lem}%%%%%%%%%%%%%%%%%%%%%%%%%%%%%%%%%%%%%%%%%%%%%%%%%%%%%%%%%%%%%%%%%%%%%%%%%%%%%%%%%%%%%%%%%%%%%%%%%%%%%%%%%%%%%
\label{lem4-2.1}
Let $H\in \bigcup _{0\leq i\leq 3}\G_{i}$, and let $X$ be a crush set of $H$.
Then the following hold.
\begin{enumerate}[{\upshape(i)}]
\item
If $H\in \G_{0}$, then $c_{1}(H-X)=c_{1}(H-X)+c_{3}(H-X)+\frac{2}{3}c_{5}(H-X)=|X|-1$ and $|X|\geq 4$.
\item
If $H\in \G^{(1)}_{1}\cup \G_{2}\cup \G_{3}\cup \G_{4}$, then $c_{1}(H-X)+c_{3}(H-X)+\frac{2}{3}c_{5}(H-X)=|X|-1$ and $|X|\geq 4$.
\item
If $H\in \G^{(2)}_{1}$, then $c_{1}(H-X)+c_{3}(H-X)+\frac{2}{3}c_{5}(H-X)=|X|-\frac{4}{3}$ and $|X|\geq 6$.
\item
If $H\in \G^{(3)}_{1}$, then $c_{1}(H-X)+c_{3}(H-X)+\frac{2}{3}c_{5}(H-X)=|X|-3$ and $|X|\geq 12$.
\end{enumerate}
\end{lem}
%%%%%%%%%%%%%%%%%%%%%%%%%%%%%%%%%%%%%%%%%%%%%%%%%%%%%%%%%%%%%%%%%%%%%%%%%%%%%%%%%%%%%%%%%%%%%%%%%%%%%%%%%%%%%%%%%%%%%%%%

%%%%%%%%%%%%%%%%%%%%%%%%%%%%%%%%%%%%%%%%%%%%%%%%%%%%%%%%%%%%%%%%%%%%%%%%%%%%%%%%%%%%%%%%%%%%%%%%%%%%%%%%%%%%%%%%%%%%%%%%
%%%%%%%%%%%%%%%%%%%%%%%%%%%%%%%%%%%%%%%%%%%%%%%%%%%%%%%%%%%%%%%%%%%%%%%%%%%%%%%%%%%%%%%%%%%%%%%%%%%%%%%%%%%%%%%%%%%%%%%%
\subsection{Hypomatchable graphs having no $\{P_{2},P_{7}\}$-factor}\label{sec4.3}
%%%%%%%%%%%%%%%%%%%%%%%%%%%%%%%%%%%%%%%%%%%%%%%%%%%%%%%%%%%%%%%%%%%%%%%%%%%%%%%%%%%%%%%%%%%%%%%%%%%%%%%%%%%%%%%%%%%%%%%%
%%%%%%%%%%%%%%%%%%%%%%%%%%%%%%%%%%%%%%%%%%%%%%%%%%%%%%%%%%%%%%%%%%%%%%%%%%%%%%%%%%%%%%%%%%%%%%%%%%%%%%%%%%%%%%%%%%%%%%%%

In this subsection, we prove the following proposition, Proposition~\ref{prop4-3.1}, which characterizes hypomatchable graphs with no $\{P_{2},P_{7}\}$-factor.
The proposition can be derived as a corollary of Proposition~\ref{prop4-4.1}, which will be proved in Subsection~\ref{sec4.4}, but we here give a proof which does not depend on Proposition~\ref{prop4-4.1} because the proof is not too long.

\begin{prop}%%%%%%%%%%%%%%%%%%%%%%%%%%%%%%%%%%%%%%%%%%%%%%%%%%%%%%%%%%%%%%%%%%%%%%%%%%%%%%%%%%%%%%%%%%%%%%%%%%%%%%%%%%%%
\label{prop4-3.1}
Let $G$ be a hypomatchable graph of order at least $7$ having no $\{P_{2},P_{7}\}$-factor.
Then $G\in \G_{0}$.
\end{prop}
%%%%%%%%%%%%%%%%%%%%%%%%%%%%%%%%%%%%%%%%%%%%%%%%%%%%%%%%%%%%%%%%%%%%%%%%%%%%%%%%%%%%%%%%%%%%%%%%%%%%%%%%%%%%%%%%%%%%%%%%
\proof
By Lemma~\ref{Thm4-1.A}, $G$ has an ear decomposition $\H=(H_{1},\ldots ,H_{m})$.
Choose $\H$ so that (H1) holds.
We use the notation introduced in Subsection~\ref{sec4.1}.

By Lemma~\ref{lem4-1.1}, $\{1\}$ is not a $7$-large set.
Hence $|V(H_{1})|\leq 5$.
Since $|V(H)|\geq 7$ by assumption, this implies $m\geq 2$.
By the definition of an ear decomposition, $H_{1}\cup H_{2}$ contains a spanning path.
Since $\{1,2\}$ is not $7$-large by Lemma~\ref{lem4-1.1}, we get $|V(H_{1})|+|V(P_{\H}(2))|\leq 5$.
Hence $|V(H_{1})|=3$.
We also have $m\geq 3$.

By Lemma~\ref{lem4-1.21}, each $H_{i}~(1\leq i\leq m)$ is a cycle of order $3$, and $G=H_{1}\cup \cdots \cup H_{m}$.
Since $|V(H)|\geq 7$ by assumption, it suffices to show that $G\in \G^{*}_{0}$.
We actually prove that for each $i~(2\leq i\leq m)$, we have $H_{1}\cup \cdots \cup H_{i}\in \G^{*}_{0}$, i.e., $H_{1}\cup \cdots \cup H_{i}\simeq K_{1}+iK_{2}$.
We proceed by induction on $i$.
We clearly have $H_{1}\cup H_{2}\simeq K_{1}+2K_{2}$.
Thus let $i\geq 3$, and assume that $H_{1}\cup \cdots \cup H_{i-1}\simeq K_{1}+(i-1)K_{2}$.
Write $V(H_{1})\cap \cdots \cap V(H_{i-1})=\{u\}$.
Suppose that $V(H_{i})\cap (V(H_{1})\cup \cdots \cup V(H_{i-1}))\not=\{u\}$.
In view of Lemma~\ref{lem4-1.00}, by relabeling $H_{1},\ldots ,H_{i-1}$ if necessary, we may assume that $V(H_{i})\cap V(H_{1})\not=\emptyset $.
Then $H_{2}\cup H_{1}\cup H_{i}$ contains a spanning path, and hence $\{1,2,i\}$ is $7$-large, which contradicts Lemma~\ref{lem4-1.1}.
Thus $V(H_{i})\cap (V(H_{1})\cup \cdots \cup V(H_{i-1}))=\{u\}$, and hence $H_{1}\cup \cdots \cup H_{i-1}\cup H_{i}\simeq K_{1}+iK_{2}$, as desired.
\qed

%%%%%%%%%%%%%%%%%%%%%%%%%%%%%%%%%%%%%%%%%%%%%%%%%%%%%%%%%%%%%%%%%%%%%%%%%%%%%%%%%%%%%%%%%%%%%%%%%%%%%%%%%%%%%%%%%%%%%%%%
%%%%%%%%%%%%%%%%%%%%%%%%%%%%%%%%%%%%%%%%%%%%%%%%%%%%%%%%%%%%%%%%%%%%%%%%%%%%%%%%%%%%%%%%%%%%%%%%%%%%%%%%%%%%%%%%%%%%%%%%
\subsection{Hypomatchable graphs having no $\{P_{2},P_{9}\}$-factor}\label{sec4.4}
%%%%%%%%%%%%%%%%%%%%%%%%%%%%%%%%%%%%%%%%%%%%%%%%%%%%%%%%%%%%%%%%%%%%%%%%%%%%%%%%%%%%%%%%%%%%%%%%%%%%%%%%%%%%%%%%%%%%%%%%
%%%%%%%%%%%%%%%%%%%%%%%%%%%%%%%%%%%%%%%%%%%%%%%%%%%%%%%%%%%%%%%%%%%%%%%%%%%%%%%%%%%%%%%%%%%%%%%%%%%%%%%%%%%%%%%%%%%%%%%%

\begin{prop}%%%%%%%%%%%%%%%%%%%%%%%%%%%%%%%%%%%%%%%%%%%%%%%%%%%%%%%%%%%%%%%%%%%%%%%%%%%%%%%%%%%%%%%%%%%%%%%%%%%%%%%%%%%%
\label{prop4-4.1}
Let $G$ be a hypomatchable graph of order at least $9$ having no $\{P_{2},P_{9}\}$-factor.
Then $G\in \G_{1}\cup \G_{2}\cup \G_{3}\cup \G_{4}$.
\end{prop}
%%%%%%%%%%%%%%%%%%%%%%%%%%%%%%%%%%%%%%%%%%%%%%%%%%%%%%%%%%%%%%%%%%%%%%%%%%%%%%%%%%%%%%%%%%%%%%%%%%%%%%%%%%%%%%%%%%%%%%%%
\proof
By Lemma~\ref{Thm4-1.A}, $G$ has an ear decomposition $\H=(H_{1},\ldots ,H_{m})$.
Choose $\H$ so that (H1) holds.

By Lemma~\ref{lem4-1.1}, $\{1\}$ is not a $9$-large set.
Hence $|V(H_{1})|\leq 7$.
Since $|V(G)|\geq 9$, this implies $m\geq 2$.
By the definition of an ear decomposition, $H_{1}\cup H_{2}$ contains a spanning path.
Since $\{1,2\}$ is not $9$-large by Lemma~\ref{lem4-1.1}, we get $|V(H_{1})|+|V(P_{\H}(2))|\leq 7$.
Hence $|V(H_{1})|=3$ or $5$.
We also have $m\geq 3$.

\medskip
\noindent
\textbf{Case 1:} $|V(H_{1})|=3$.

By Lemma~\ref{lem4-1.21}, each $H_{i}~(1\leq i\leq 3)$ is a cycle of order $3$, and $G=H_{1}\cup \cdots \cup H_{m}$.
We show that $G\in \G_{1}$.
We actually prove that for each $i~(2\leq i\leq m)$, we have $H_{1}\cup \cdots \cup H_{i}\in \G^{*}_{1}$, i.e., $H_{1}\cup \cdots \cup H_{i}\simeq A_{1}(s_{1},s_{2},s_{3})$ for some $s_{1},s_{2},s_{3}$ with $s_{1}+s_{2}+s_{3}=i-1$.
We proceed by induction on $i$.
Note that $H_{1}\cup H_{2}\simeq A_{1}(1,0,0)$.
Thus let $i\geq 3$, and assume that $H_{1}\cup \cdots \cup H_{i-1}\simeq A_{1}(s'_{1},s'_{2},s'_{3})$ with $s'_{1}+s'_{2}+s'_{3}=i-2$.
If only one of $s'_{1}$, $s'_{2}$ and $s'_{3}$ is nonzero, i.e., $H_{1}\cup \cdots \cup H_{i-1}\simeq A_{1}(i-2,0,0)$, then $H_{1}\cup \cdots \cup H_{i-1}\cup H_{i}\simeq A_{1}(i-1,0,0)$ or $A_{1}(i-2,1,0)$.
Thus we may assume that at least two of $s'_{1}$, $s'_{2}$ and $s'_{3}$ are nonzero.
In view of Lemma~\ref{lem4-1.00}, by relabeling $H_{1},\ldots ,H_{i-1}$ if necessary, we may assume that $H_{1}$ intersects with all of $H_{2},\ldots ,H_{i-1}$.
Write $H_{1}=w_{1}w_{2}w_{3}w_{1}$.
We may assume that $s'_{h}=|\{j\mid 2\leq j\leq i-1,V(H_{j})\cap V(H_{1})=\{w_{h}\}\}|$ for each $h=1,2,3$.
Suppose that there exists $j~(2\leq j\leq i-1)$ such that $V(H_{i})\cap V(P_{\H}(j))\not=\emptyset $.
Since at least two of $s'_{1},s'_{2},s'_{3}$ are nonzero, there exists $j'~(2\leq j'\leq i-1)$ with $j'\not=j$ such that $V(H_{j'})\cap V(H_{1})\not=V(H_{j})\cap V(H_{1})$.
Then $H_{j'}\cup H_{1}\cup H_{j}\cup H_{i}$ contains a spanning path, and hence $\{1,j,j',i\}$ is $9$-large, which contradicts Lemma~\ref{lem4-1.1}.
Consequently $V(H_{i})\cap V(P_{\H}(j))=\emptyset $ for every $j~(2\leq j\leq i-1)$, which implies $V(H_{1})\cap (\bigcup _{1\leq j\leq i-1}V(H_{j}))=V(H_{i})\cap V(H_{1})$.
We may assume $V(H_{i})\cap V(H_{1})=\{w_{1}\}$.
Thus $H_{1}\cup \cdots \cup H_{i-1}\cup H_{i}\simeq A_{1}(s'_{1}+1,s'_{2},s'_{3})$, as desired.

\medskip
\noindent
\textbf{Case 2:} $|V(H_{1})|=5$.

We first prove two claims.

\begin{claim}%%%%%%%%%%%%%%%%%%%%%%%%%%%%%%%%%%%%%%%%%%%%%%%%%%%%%%%%%%%%%%%%%%%%%%%%%%%%%%%%%%%%%%%%%%%%%%%%%%%%%%%%%%%
\label{cl4-4.1.1}
For each $i~(2\leq i\leq m)$, $|V(P_{\H}(i))|=2$ and $N_{G}(V(P_{\H}(i)))\cap (\bigcup _{2\leq j\leq i-1}V(P_{\H}(j)))=\emptyset $.
\end{claim}
%%%%%%%%%%%%%%%%%%%%%%%%%%%%%%%%%%%%%%%%%%%%%%%%%%%%%%%%%%%%%%%%%%%%%%%%%%%%%%%%%%%%%%%%%%%%%%%%%%%%%%%%%%%%%%%%%%%%%%%%
\proof
We proceed by induction on $i$.
Let $i\geq 2$, and assume that for each $i'$ with $2\leq i'\leq i-1$, we have $|V(P_{\H}(i'))|=2$ and $N_{G}(V(P_{\H}(i')))\cap (\bigcup _{2\leq j\leq i'-1}V(P_{\H}(j)))=\emptyset $ (this includes the case where $i=2$).
It follows from (E4) that for each $i'~(2\leq i'\leq i-1)$ and for each $v\in V(P_{\H}(i'))$, $H_{1}\cup H_{i'}$ contains a spanning path having $v$ as one of its endvertices.
Let $U$ be the set of the endvertices of $P_{\H}(i)$.
Suppose that $N_{G}(U)\cap (\bigcup _{2\leq j\leq i-1}V(P_{\H}(j)))\not=\emptyset $, and take $v\in N_{G}(U)\cap (\bigcup _{2\leq j\leq i-1}V(P_{\H}(j)))$.
Let $i'$ denote the index such that $v\in V(P_{\H}(i'))$.
Then since $H_{1}\cup H_{i'}$ contains a spanning path having endvertex $v$, $G[V(H_{1})\cup V(P_{\H}(i'))\cup V(P_{\H}(i))]$ contains a spanning path.
Since $|V(H_{1})|+|V(P_{\H}(i'))|+|V(P_{\H}(i))|=7+|V(P_{\H}(i))|\geq 9$, this contradicts Lemma~\ref{lem4-1.1}.
Thus $N_{G}(U)\cap (\bigcup _{2\leq j\leq i-1}V(P_{\H}(j)))=\emptyset $.
It now follows from (E4) that $H_{1}\cup H_{i}$ contains a spanning path.
Hence by Lemma~\ref{lem4-1.1}, $|V(P_{\H}(i))|\leq 7-|V(H_{1})|=2$.
This implies $U=V(P_{\H}(i))$, and thus the claim is proved.
\qed

\begin{claim}%%%%%%%%%%%%%%%%%%%%%%%%%%%%%%%%%%%%%%%%%%%%%%%%%%%%%%%%%%%%%%%%%%%%%%%%%%%%%%%%%%%%%%%%%%%%%%%%%%%%%%%%%%%
\label{cl4-4.1.2}
If $N_{G}(\bigcup _{2\leq i\leq m}V(P_{\H}(i)))\cap V(H_{1})$ contains two vertices $w,w'$ which are consecutive on $H_{1}$, then $|N_{G}(w)\cap (\bigcup _{2\leq i\leq m}V(P_{\H}(i)))|=|N_{G}(w')\cap (\bigcup _{2\leq i\leq m}V(P_{\H}(i)))|=1$ and $N_{G}(w)\cap (\bigcup _{2\leq i\leq m}V(P_{\H}(i)))=N_{G}(w')\cap (\bigcup _{2\leq i\leq m}V(P_{\H}(i)))$.
\end{claim}
%%%%%%%%%%%%%%%%%%%%%%%%%%%%%%%%%%%%%%%%%%%%%%%%%%%%%%%%%%%%%%%%%%%%%%%%%%%%%%%%%%%%%%%%%%%%%%%%%%%%%%%%%%%%%%%%%%%%%%%%
\proof
Suppose that $|N_{G}(w)\cap (\bigcup _{2\leq i\leq m}V(P_{\H}(i)))|\geq 2$ or $|N_{G}(w')\cap (\bigcup _{2\leq i\leq m}V(P_{\H}(i)))|\geq 2$ or $N_{G}(w)\cap (\bigcup _{2\leq i\leq m}V(P_{\H}(i)))\not=N_{G}(w')\cap (\bigcup _{2\leq i\leq m}V(P_{\H}(i)))$.
Then we can take $v\in N_{G}(w)\cap (\bigcup _{2\leq i\leq m}V(P_{\H}(i)))$ and $v'\in N_{G}(w')\cap (\bigcup _{2\leq i\leq m}V(P_{\H}(i)))$ so that $v\not=v'$.
Let $i$ and $i'$ be the indices such that $v\in V(P_{\H}(i))$ and $v'\in V(P_{\H}(i'))$.
By Claim~\ref{cl4-4.1.1}, $|V(P_{\H}(i))|=|V(P_{\H}(i'))|=2$.
Hence by Lemma~\ref{lem4-1.20}, $i\not=i'$.
Note that $G[V(H_{1})\cup V(P_{\H}(i))\cup V(P_{\H}(i'))]$ contains a spanning path.
Since $|V(H_{1})|+|V(P_{\H}(i))|+|V(P_{\H}(i'))|=9$, this contradicts Lemma~\ref{lem4-1.1}.
\qed

We return to the proof of Proposition~\ref{prop4-4.1}.
Write $H_{1}=w_{1}w_{2}w_{3}w_{4}w_{5}w_{1}$.
We first consider the case where $N_{G}(\bigcup _{1\leq i\leq m}V(P_{\H}(i)))\cap V(H_{1})$ contains two vertices $w,w'$ which are consecutive on $H_{1}$.
We may assume $w=w_{3}$ and $w'=w_{4}$.
By Claim~\ref{cl4-4.1.2}, there exists $b\in \bigcup _{1\leq i\leq m}V(P_{\H}(i))$ such that $N_{G}(w_{3})\cap (\bigcup _{1\leq i\leq m}V(P_{\H}(i)))=N_{G}(w_{4})\cap (\bigcup _{1\leq i\leq m}V(P_{\H}(i)))=\{b\}$.
Note that Claim~\ref{cl4-4.1.1} in particular implies that for any permutation $i_{2},\ldots ,i_{m}$ of $2,\ldots ,m$, $(H_{1},H_{i_{2}},\ldots ,H_{i_{m}})$ is an ear decomposition.
Thus we may assume $b\in V(P_{\H}(2))$.
Write $P_{\H}(2)=bb'$.
By Claim~\ref{cl4-4.1.2} and (E4), $N_{G}(b')\cap V(H_{1})=\{w_{1}\}$, $\{w_{3},w_{4}\}\subseteq N_{G}(b)\cap V(H_{1})\subseteq \{w_{1},w_{3},w_{4}\}$, and $N_{G}(v)\cap V(H_{1})=\{w_{1}\}$ for all $v\in \bigcup _{3\leq i\leq m}V(P_{\H}(i))$.
Consequently $A'_{4}(m-2)\subseteq G$.
By Lemma~\ref{lem4-1.1}, $\{1,2,3\}$ is not $9$-large.
Hence $G[V(H_{1})\cup V(P_{\H}(2))]$ does not contain a spanning path with endvertex $w_{1}$.
This implies $w_{2}w_{4},w_{2}w_{5},w_{3}w_{5}\not\in E(G)$, and hence it follows from Claim~\ref{cl4-4.1.1} that $G\subseteq A''_{4}(m-2)$.
Therefore $G\in \G_{4}$.

We now consider the case where $N_{G}(\bigcup _{2\leq i\leq m}V(P_{\H}(i)))\cap V(H_{1})$ does not contain two consecutive vertices.
%In this case, $N_{G}(\bigcup _{2\leq i\leq m}V(P_{\H}(i)))\cap V(H_{1})$ consists of two nonconsecutive vertices.
In this case, $|N_{G}(\bigcup _{2\leq i\leq m}V(P_{\H}(i)))\cap V(H_{1})|\leq 2$.
% consists of at most two vertices, and if it consists of exactly two vertices, then the vertices are not consecutive on $H_{1}$.
We may assume $N_{G}(\bigcup _{2\leq i\leq m}V(P_{\H}(i)))\cap V(H_{1})\subseteq \{w_{1},w_{3}\}$.
Let $s_{h}=|\{i\mid 2\leq i\leq m,N_{G}(V(P_{\H}(i)))\cap V(H_{1})=\{w_{h}\}\}|$ for $h=1,3$, and $s_{2}=|\{i\mid 2\leq i\leq m,N_{G}(V(P_{\H}(i)))\cap V(H_{1})=\{w_{1},w_{3}\}\}|$.
Since $m\geq 3$, $s_{1}+s_{2}+s_{3}\geq 2$.
If $s_{2}=0$ and $s_{1}$ or $s_{3}$ (say $s_{3}$) is zero, then it follows from Claim~\ref{cl4-4.1.1} that $A'_{3}(s_{1})\subseteq G\subseteq A''_{3}(s_{1})$, and hence $G\in \G_{3}$.
Thus we may assume that we have $s_{2}\not=0$, or $s_{1}\not=0$ and $s_{3}\not=0$.
Since $s_{1}+s_{2}+s_{3}\geq 2$, it follows from Lemma~\ref{lem4-1.1} that $G[V(H_{1})]$ does not contain a spanning path connecting $w_{1}$ and $w_{3}$.
Hence $w_{2}w_{4},w_{2}w_{5}\not\in E(G)$, which together with Claim~\ref{cl4-4.1.1} implies that $A'_{2}(s_{1},s_{2}+1,s_{3})\subseteq G\subseteq A''_{2}(s_{1},s_{2}+1,s_{3})$.
Therefore $G\in \G_{2}$.

This completes the proof of Proposition~\ref{prop4-4.1}.
\qed

%%%%%%%%%%%%%%%%%%%%%%%%%%%%%%%%%%%%%%%%%%%%%%%%%%%%%%%%%%%%%%%%%%%%%%%%%%%%%%%%%%%%%%%%%%%%%%%%%%%%%%%%%%%%%%%%%%%%%%%%
%%%%%%%%%%%%%%%%%%%%%%%%%%%%%%%%%%%%%%%%%%%%%%%%%%%%%%%%%%%%%%%%%%%%%%%%%%%%%%%%%%%%%%%%%%%%%%%%%%%%%%%%%%%%%%%%%%%%%%%%
\subsection{Alternating paths}\label{sec4.5}
%%%%%%%%%%%%%%%%%%%%%%%%%%%%%%%%%%%%%%%%%%%%%%%%%%%%%%%%%%%%%%%%%%%%%%%%%%%%%%%%%%%%%%%%%%%%%%%%%%%%%%%%%%%%%%%%%%%%%%%%
%%%%%%%%%%%%%%%%%%%%%%%%%%%%%%%%%%%%%%%%%%%%%%%%%%%%%%%%%%%%%%%%%%%%%%%%%%%%%%%%%%%%%%%%%%%%%%%%%%%%%%%%%%%%%%%%%%%%%%%%

In this appendant subsection, we prove two lemmas about hypomatchable graphs, which we use in the proof of Theorem~\ref{mainthm2}.
Throughout this subsection, we let $G$ denote a hypomatchable graph, let $v\in V(G)$, and let $M$ be a perfect matching of $G-v$.
A path $v_{1}v_{2}\cdots v_{l}$ with $v_{1}=v$ is called an {\it alternating path} if $v_{2i}v_{2i+1}\in M$ for each $i$ with $1\leq i\leq \frac{l-1}{2}$.

\begin{lem}%%%%%%%%%%%%%%%%%%%%%%%%%%%%%%%%%%%%%%%%%%%%%%%%%%%%%%%%%%%%%%%%%%%%%%%%%%%%%%%%%%%%%%%%%%%%%%%%%%%%%%%%%%%%%
\label{lem4-5.1}
For each $w\in V(G)$, $G$ contains an alternating path $Q$ of odd order connecting $v$ and $w$ such that $M-E(Q)$ is a perfect matching of $G-V(Q)$.
\end{lem}
%%%%%%%%%%%%%%%%%%%%%%%%%%%%%%%%%%%%%%%%%%%%%%%%%%%%%%%%%%%%%%%%%%%%%%%%%%%%%%%%%%%%%%%%%%%%%%%%%%%%%%%%%%%%%%%%%%%%%%%%
\proof
If $w=v$, then it suffices simply to let $Q=v$.
Thus we may assume $w\not=v$.
Let $M'$ be a perfect matching of $G-w$, and let $H$ denote the subgraph induced by the symmetric difference of $M$ and $M'$.
Then $d_{H}(v)=d_{H}(w)=1$, and $d_{H}(x)=2$ for all $x\in V(H)-\{v,w\}$.
This implies that the component $Q$ of $H$ containing $v$ is an alternating path connecting $v$ and $w$.
Since the edge of $Q$ incident with $v$ does not belong to $M$ and the edge of $Q$ incident with $w$ belongs to $M$, $Q$ has odd order, and $M-E(Q)$ is a perfect matching of $G-V(Q)$.
\qed

\begin{lem}%%%%%%%%%%%%%%%%%%%%%%%%%%%%%%%%%%%%%%%%%%%%%%%%%%%%%%%%%%%%%%%%%%%%%%%%%%%%%%%%%%%%%%%%%%%%%%%%%%%%%%%%%%%%%
\label{lem4-5.2}
Suppose that $|V(G)|\geq 5$ and, in the case where $G$ is isomorphic to $K_{1}+sK_{2}$ for some $s\geq 2$, suppose further that $v$ is not the unique cutvertex of $G$.
Then $G$ contains an alternating path $Q$ of odd order having $v$ as one of its endvertices such that $|V(Q)|\geq 5$ and $M-E(Q)$ is a perfect matching of $G-V(Q)$.
\end{lem}
%%%%%%%%%%%%%%%%%%%%%%%%%%%%%%%%%%%%%%%%%%%%%%%%%%%%%%%%%%%%%%%%%%%%%%%%%%%%%%%%%%%%%%%%%%%%%%%%%%%%%%%%%%%%%%%%%%%%%%%%
\proof
If $vu\in E(G)$ for all $u\in V(G)-\{v\}$, then the assumption of the lemma implies that $G$ contains an edge $xy$ joining endvertices of two distinct edges $xx',yy'$ in $M$, and hence $vx'xyy'$ is a path with the desired properties.
Thus we may assume that there exists $u\in V(G)-\{v\}$ such that $vu\not\in E(G)$.
Let $uw\in M$.
By Lemma~\ref{lem4-5.1}, $G$ contains an alternating path $Q$ of odd order connecting $v$ and $w$ such that $M-E(Q)$ is a perfect matching of $G-V(Q)$.
Since $Q$ is an alternating path of odd order and $vu\not\in E(G)$, we get $|V(Q)|\geq 5$, as desired.
\qed

%%%%%%%%%%%%%%%%%%%%%%%%%%%%%%%%%%%%%%%%%%%%%%%%%%%%%%%%%%%%%%%%%%%%%%%%%%%%%%%%%%%%%%%%%%%%%%%%%%%%%%%%%%%%%%%%%%%%%%%%
%%%%%%%%%%%%%%%%%%%%%%%%%%%%%%%%%%%%%%%%%%%%%%%%%%%%%%%%%%%%%%%%%%%%%%%%%%%%%%%%%%%%%%%%%%%%%%%%%%%%%%%%%%%%%%%%%%%%%%%%
%%%%%%%%%%%%%%%%%%%%%%%%%%%%%%%%%%%%%%%%%%%%%%%%%%%%%%%%%%%%%%%%%%%%%%%%%%%%%%%%%%%%%%%%%%%%%%%%%%%%%%%%%%%%%%%%%%%%%%%%
\section{Proof of main theorems}\label{sec5}
%%%%%%%%%%%%%%%%%%%%%%%%%%%%%%%%%%%%%%%%%%%%%%%%%%%%%%%%%%%%%%%%%%%%%%%%%%%%%%%%%%%%%%%%%%%%%%%%%%%%%%%%%%%%%%%%%%%%%%%%
%%%%%%%%%%%%%%%%%%%%%%%%%%%%%%%%%%%%%%%%%%%%%%%%%%%%%%%%%%%%%%%%%%%%%%%%%%%%%%%%%%%%%%%%%%%%%%%%%%%%%%%%%%%%%%%%%%%%%%%%
%%%%%%%%%%%%%%%%%%%%%%%%%%%%%%%%%%%%%%%%%%%%%%%%%%%%%%%%%%%%%%%%%%%%%%%%%%%%%%%%%%%%%%%%%%%%%%%%%%%%%%%%%%%%%%%%%%%%%%%%

For a graph $H$, we let $\C_{\rm odd}(H)$ denote the set of those components of $H$ having odd order, and set $c_{\rm odd}(H)=|\C_{\rm odd}(H)|$.

Recall that Tutte's $1$-factor theorem says that if a graph $G$ of even order has no perfect matching, then there exists $S\subseteq V(G)$ such that $c_{\rm odd}(G-S)\geq |S|+2$.
In this section, we often choose a set $S$ of vertices of a given graph $G$ so that
\begin{enumerate}[{\bf (S1)}]
\item
$c_{\rm odd}(G-S)-|S|$ is as large as possible, and
\item
subject to (S1), $|S|$ is as large as possible.
\end{enumerate}
Note that $c_{\rm odd}(G-S)-|S|\geq c_{\rm odd}(G)-|\emptyset |\geq 0$ (it is possible that $S=\emptyset $, but our argument in this section works even if $S=\emptyset $).

We first give a fundamental lemma.

\begin{lem}%%%%%%%%%%%%%%%%%%%%%%%%%%%%%%%%%%%%%%%%%%%%%%%%%%%%%%%%%%%%%%%%%%%%%%%%%%%%%%%%%%%%%%%%%%%%%%%%%%%%%%%%%%%%%
\label{lem5.1}
Let $G$ be a graph, and let $S$ be a subset of $V(G)$ satisfying (S1) and (S2).
Then the following hold.
\begin{enumerate}[{\upshape(i)}]
\item
We have $\C(G-S)=\C_{\rm odd}(G-S)$.
\item
For each $C\in \C_{\rm odd}(G-S)$, $C$ is hypomatchable.
\item
Let $H$ be the bipartite graph $H$ with bipartition $(S,\C_{\rm odd}(G-S))$ defined by letting $uC\in E(H)~(u\in S,C\in \C_{\rm odd}(G-S))$ if and only if $N_{G}(u)\cap V(C)\not=\emptyset $.
Then for every $X\subseteq S$, $|N_{H}(X)|\geq |X|$.
\end{enumerate}
\end{lem}
%%%%%%%%%%%%%%%%%%%%%%%%%%%%%%%%%%%%%%%%%%%%%%%%%%%%%%%%%%%%%%%%%%%%%%%%%%%%%%%%%%%%%%%%%%%%%%%%%%%%%%%%%%%%%%%%%%%%%%%%
\proof
\begin{enumerate}[{\upshape(i)}]
\item
Suppose that there exists $C\in \C(G-S)$ such that $|V(C)|$ is even, and take $v\in V(C)$.
Then $c_{\rm odd}(C-v)\geq 1$.
Let $S_{1}=S\cup \{v\}$.
Then $c_{\rm odd}(G-S_{1})-|S_{1}|=(c_{\rm odd}(G-S)+c_{\rm odd}(C-v))-(|S|+1)\geq c_{\rm odd}(G-S)-|S|$ and $|S_{1}|>|S|$, which contradicts (S1) or (S2).
\item
Suppose that $C$ is not hypomatchable.
Then there exists $v\in V(C)$ such that $C-v$ has no perfect matching.
Applying Tutte's $1$-factor theorem to $C-v$, we see that there exists $S''\subseteq V(C)$ with $v\in S''$ such that $c_{\rm odd}(C-S'')\geq |S''|+1$.
Let $S_{2}=S\cup S''$.
Then $c_{\rm odd}(G-S_{2})-|S_{2}|=(c_{\rm odd}(G-S)-1+c_{\rm odd}(C-S''))-(|S|+|S''|)\geq c_{\rm odd}(G-S)-|S|$ and $|S_{2}|=|S|+|S''|>|S|$, which contradicts (S1) or (S2).
\item
Suppose that there exists $X\subseteq S$ such that $|N_{H}(X)|<|X|$.
Set $S_{3}=S-X$.
Then every component in $\C_{\rm odd}(G-S)-N_{H}(X)$ belongs to $\C_{\rm odd}(G-S_{3})$.
Hence
\begin{align*}
c_{\rm odd}(G-S_{3})-|S_{3}| &\geq (c_{\rm odd}(G-S)-|N_{H}(X)|)-|S_{3}|\\
&> c_{\rm odd}(G-S)-|X|-|S_{3}|\\
&= c_{\rm odd}(G-S)-|S|,
\end{align*}
which contradicts (S1).
\qed
\end{enumerate}

%%%%%%%%%%%%%%%%%%%%%%%%%%%%%%%%%%%%%%%%%%%%%%%%%%%%%%%%%%%%%%%%%%%%%%%%%%%%%%%%%%%%%%%%%%%%%%%%%%%%%%%%%%%%%%%%%%%%%%%%
%%%%%%%%%%%%%%%%%%%%%%%%%%%%%%%%%%%%%%%%%%%%%%%%%%%%%%%%%%%%%%%%%%%%%%%%%%%%%%%%%%%%%%%%%%%%%%%%%%%%%%%%%%%%%%%%%%%%%%%%
\subsection{Proof of Theorem~\ref{mainthm1}}\label{sec5.1}
%%%%%%%%%%%%%%%%%%%%%%%%%%%%%%%%%%%%%%%%%%%%%%%%%%%%%%%%%%%%%%%%%%%%%%%%%%%%%%%%%%%%%%%%%%%%%%%%%%%%%%%%%%%%%%%%%%%%%%%%
%%%%%%%%%%%%%%%%%%%%%%%%%%%%%%%%%%%%%%%%%%%%%%%%%%%%%%%%%%%%%%%%%%%%%%%%%%%%%%%%%%%%%%%%%%%%%%%%%%%%%%%%%%%%%%%%%%%%%%%%

For a graph $H$, we let $\C'(H)$ denote the set of those components $C\in \C_{\rm odd}(H)$ such that $|V(C)|\geq 3$ and $C$ is a hypomatchable graph having no $\{P_{2},P_{7}\}$-factor, and set $c'(H)=|\C'(H)|$.

We first give a sufficient condition for the existence of a $\{P_{2},P_{7}\}$-factor in terms of $c_{1}$ and $c'$.

\begin{thm}%%%%%%%%%%%%%%%%%%%%%%%%%%%%%%%%%%%%%%%%%%%%%%%%%%%%%%%%%%%%%%%%%%%%%%%%%%%%%%%%%%%%%%%%%%%%%%%%%%%%%%%%%%%%%
\label{thm5-1.1}
Let $G$ be a graph.
If $c_{1}(G-X)+\frac{1}{2}c'(G-X)\leq |X|$ for all $X\subseteq V(G)$, then $G$ has a $\{P_{2},P_{7}\}$-factor.
\end{thm}
%%%%%%%%%%%%%%%%%%%%%%%%%%%%%%%%%%%%%%%%%%%%%%%%%%%%%%%%%%%%%%%%%%%%%%%%%%%%%%%%%%%%%%%%%%%%%%%%%%%%%%%%%%%%%%%%%%%%%%%%
\proof
Choose $S\subseteq V(G)$ so that (S1) and (S2) hold.

Set $T=\C_{\rm odd}(G-S)~(=\C(G-S))$, $T_{1}=\C_{1}(G-S)$ and $T_{2}=\C'(G-S)$.
Then $T_{1}\cap T_{2}=\emptyset $ and $T_{1}\cup T_{2}\subseteq T$.
We construct a bipartite graph $H$ with bipartition $(S,T)$ by letting $uC\in E(H)~(u\in S,C\in T)$ if and only if $N_{G}(u)\cap V(C)\not=\emptyset $.

\begin{claim}%%%%%%%%%%%%%%%%%%%%%%%%%%%%%%%%%%%%%%%%%%%%%%%%%%%%%%%%%%%%%%%%%%%%%%%%%%%%%%%%%%%%%%%%%%%%%%%%%%%%%%%%%%%
\label{cl5-1.1.1}
For every $Y\subseteq T_{1}\cup T_{2}$, $|N_{H}(Y)|\geq |Y\cap T_{1}|+\frac{1}{2}|Y\cap T_{2}|$.
\end{claim}
%%%%%%%%%%%%%%%%%%%%%%%%%%%%%%%%%%%%%%%%%%%%%%%%%%%%%%%%%%%%%%%%%%%%%%%%%%%%%%%%%%%%%%%%%%%%%%%%%%%%%%%%%%%%%%%%%%%%%%%%
\proof
Suppose that there exists $Y\subseteq T_{1}\cup T_{2}$ such that $|N_{H}(Y)|<|Y\cap T_{1}|+\frac{1}{2}|Y\cap T_{2}|$.
Set $X'=N_{H}(Y)$.
Then each element of $Y\cap T_{1}$ belongs to $\C_{1}(G-X')$, and each element of $Y\cap T_{2}$ belongs to $\C'(G-X')$.
Hence $|Y\cap T_{1}|\leq c_{1}(G-X')$ and $|Y\cap T_{2}|\leq c'(G-X')$.
Consequently $|X'|=|N_{H}(Y)|<|Y\cap T_{1}|+\frac{1}{2}|Y\cap T_{2}|\leq c_{1}(G-X')+\frac{1}{2}c'(G-X')$, which contradicts the assumption of the theorem.
\qed

Now we apply Proposition~\ref{prop3.1} with $G$ and $L$ replaced by $H$ and $\emptyset $, respectively.
Then by Lemma~\ref{lem5.1}(iii) and Claim~\ref{cl5-1.1.1}, $H$ has a subgraph $F$ with $V(F)\supseteq S\cup T_{1}\cup T_{2}$ such that each $A\in \C(F)$ is a path satisfying one of (I) and (II) in Proposition~\ref{prop3.1}.
For $A\in \C(F)$, let $U_{A}=V(A)\cap S$ and $\L_{A}=V(A)\cap T$, and let $G_{A}=G[U_{A}\cup (\bigcup _{C\in \L_{A}}V(C))]$.

\begin{claim}%%%%%%%%%%%%%%%%%%%%%%%%%%%%%%%%%%%%%%%%%%%%%%%%%%%%%%%%%%%%%%%%%%%%%%%%%%%%%%%%%%%%%%%%%%%%%%%%%%%%%%%%%%%
\label{cl5-1.1.2}
For each $A\in \C(F)$, $G_{A}$ has a $\{P_{2},P_{7}\}$-factor.
\end{claim}
%%%%%%%%%%%%%%%%%%%%%%%%%%%%%%%%%%%%%%%%%%%%%%%%%%%%%%%%%%%%%%%%%%%%%%%%%%%%%%%%%%%%%%%%%%%%%%%%%%%%%%%%%%%%%%%%%%%%%%%%
\proof
We first assume that $A$ satisfies (I).
Then $|U_{A}|=|\L_{A}|=1$.
Write $U_{A}=\{u\}$ and $\L_{A}=\{D\}$, and let $v\in V(D)$ be a vertex with $uv\in E(G)$.
Since $D$ is hypomatchable by Lemma~\ref{lem5.1}(ii), $D-v$ has a perfect matching $M$.
Hence $M\cup \{uv\}$ is a perfect matching of $G_{A}$.
In particular, $G_{A}$ has a $\{P_{2},P_{7}\}$-factor.

Next we assume that $A$ satisfies (II).
Note that $|V(A)|$ is odd and $|V(A)|\geq 3$.
Write $A=D_{1}u_{1}D_{2}u_{2}\cdots D_{l}u_{l}D_{l+1}~(u_{i}\in U_{A},D_{i}\in \L_{A})$.
Let $v_{i}\in N_{G}(u_{i})\cap V(D_{i})$ for $1\leq i\leq l$, and let $v_{l+1}\in N_{G}(u_{l})\cap V(D_{l+1})$.
Since $A$ satisfies (II), $|V(D_{1})-\{v_{1}\}|\geq 2$, $|V(D_{l+1})-\{v_{l+1}\}|\geq 2$ and $V(D_{i})=\{v_{i}\}~(2\leq i\leq l)$.
Fix $i\in \{1,l+1\}$.
Since $D_{i}$ is hypomatchable by the definition of $T_{2}$, $D_{i}-v_{i}$ has a perfect matching $M_{i}$.
Since $|V(D_{i})|\geq 3$, $v_{i}$ is adjacent to a vertex $u'_{i}\in V(D_{i})$.
Let $v'_{i}\in V(D_{i})$ be the vertex with $u'_{i}v'_{i}\in M_{i}$.
Then $P=v'_{1}u'_{1}v_{1}u_{1}v_{2}u_{2}\cdots v_{l}u_{l}v_{l+1}u'_{l+1}v'_{l+1}$ is a path of order at least $7$.
Since $M_{i}-\{u'_{i}v'_{i}\}$ is a matching for each $i\in \{1,l+1\}$, $F_{A}=P\cup (M_{1}-\{u'_{1}v'_{1}\})\cup (M_{l+1}-\{u'_{l+1}v'_{l+1}\})$ is a path-factor of $G_{A}$ with $\C_{3}(F_{A})=\C_{5}(F_{A})=\emptyset $.
By Fact~\ref{fact1}, $G_{A}$ has a $\{P_{2},P_{7}\}$-factor.
\qed

By Lemma~\ref{lem5.1}(i)(ii), each component in $\C(G-S)-\C_{1}(G-S)-\C'(G-S)$ has a $\{P_{2},P_{7}\}$-factor.
This together with Claim~\ref{cl5-1.1.2} implies that $G$ has a $\{P_{2},P_{7}\}$-factor.

This completes the proof of Theorem~\ref{thm5-1.1}.
\qed

\medbreak\noindent\textit{Proof of Theorem~\ref{mainthm1}.}\quad
Let $G$ be as in Theorem~\ref{mainthm1}.
Suppose that $G$ has no $\{P_{2},P_{7}\}$-factor.
Then by Theorem~\ref{thm5-1.1}, there exists $X\subseteq V(G)$ such that $c_{1}(G-X)+\frac{1}{2}c'(G-X)>|X|$.
Write $\C'(G-X)-(\C_{3}(G-X)\cup \C_{5}(G-X))=\{D_{1},\ldots ,D_{q}\}$.
For each $i~(1\leq i\leq q)$, since $D_{i}$ is a hypomatchable graph of order at least $7$ with no $\{P_{2},P_{7}\}$-factor, it follows from Proposition~\ref{prop4-3.1} that $D_{i}\in \G_{0}$.
For each $i~(1\leq i\leq q)$, let $X_{i}$ be a crush set of $D_{i}$.
By Lemma~\ref{lem4-2.1}, $c_{1}(D_{i}-X_{i})=|X_{i}|-1$ and $|X_{i}|\geq 4$, and hence $c_{1}(D_{i}-X_{i})\geq \frac{3}{4}|X_{i}|$.
Let $X_{0}=X\cup (\bigcup _{1\leq i\leq q}X_{i})$.

Then $c_{1}(G-X_{0})=c_{1}(G-X)+\sum _{1\leq i\leq q}c_{1}(D_{i}-X_{i})\geq c_{1}(G-X)+\frac{3}{4}\sum _{1\leq i\leq q}|X_{i}|$.
Consequently
\begin{align*}
c_{1}(G-X_{0})&-\frac{2}{3}c_{1}(G-X)-\frac{1}{3}q-\frac{2}{3}\sum _{1\leq i\leq q}|X_{i}|\\
&\geq c_{1}(G-X)+\frac{3}{4}\sum _{1\leq i\leq q}|X_{i}|-\frac{2}{3}c_{1}(G-X)-\frac{1}{3}q-\frac{2}{3}\sum _{1\leq i\leq q}|X_{i}|\\
&= \frac{1}{3}c_{1}(G-X)+\sum _{1\leq i\leq q}\left(\frac{3}{4}|X_{i}|-\frac{1}{3}-\frac{2}{3}|X_{i}|\right)\\
&= \frac{1}{3}c_{1}(G-X)+\sum _{1\leq i\leq q}\left(\frac{1}{12}|X_{i}|-\frac{1}{3}\right)\\
&\geq 0,
\end{align*}
and hence
$$
\frac{2}{3}c_{1}(G-X)+\frac{1}{3}q+\frac{2}{3}\sum _{1\leq i\leq q}|X_{i}|\leq c_{1}(G-X_{0}).
$$
This leads to
\begin{align*}
\frac{2}{3}|X_{0}| &= \frac{2}{3}\left(|X|+\sum _{1\leq i\leq q}|X_{i}|\right)\\
&< \frac{2}{3}\left(c_{1}(G-X)+\frac{1}{2}c'(G-X)+\sum _{1\leq i\leq q}|X_{i}|\right)\\
&= \frac{2}{3}c_{1}(G-X)+\frac{1}{3}|\C'(G-X)\cap \C_{3}(G-X)|+\frac{1}{3}|\C'(G-X)\cap \C_{5}(G-X)|\\
&\quad +\frac{1}{3}|\C'(G-X)-(\C_{3}(G-X)\cup \C_{5}(G-X))|+\frac{2}{3}\sum _{1\leq i\leq q}|X_{i}|\\
&\leq \frac{2}{3}c_{1}(G-X)+\frac{1}{3}c_{3}(G-X)+\frac{1}{3}c_{5}(G-X)+\frac{1}{3}q+\frac{2}{3}\sum _{1\leq i\leq q}|X_{i}|\\
&= \frac{2}{3}c_{1}(G-X)+\frac{1}{3}c_{3}(G-X_{0})+\frac{1}{3}c_{5}(G-X_{0})+\frac{1}{3}q+\frac{2}{3}\sum _{1\leq i\leq q}|X_{i}|\\
&\leq c_{1}(G-X_{0})+\frac{1}{3}c_{3}(G-X_{0})+\frac{1}{3}c_{5}(G-X_{0}),
\end{align*}
which contradicts the assumption of the theorem.

This completes the proof of Theorem~\ref{mainthm1}.
\qed

%%%%%%%%%%%%%%%%%%%%%%%%%%%%%%%%%%%%%%%%%%%%%%%%%%%%%%%%%%%%%%%%%%%%%%%%%%%%%%%%%%%%%%%%%%%%%%%%%%%%%%%%%%%%%%%%%%%%%%%%
%%%%%%%%%%%%%%%%%%%%%%%%%%%%%%%%%%%%%%%%%%%%%%%%%%%%%%%%%%%%%%%%%%%%%%%%%%%%%%%%%%%%%%%%%%%%%%%%%%%%%%%%%%%%%%%%%%%%%%%%
\subsection{Proof of Theorem~\ref{mainthm2}}\label{sec5.2}
%%%%%%%%%%%%%%%%%%%%%%%%%%%%%%%%%%%%%%%%%%%%%%%%%%%%%%%%%%%%%%%%%%%%%%%%%%%%%%%%%%%%%%%%%%%%%%%%%%%%%%%%%%%%%%%%%%%%%%%%
%%%%%%%%%%%%%%%%%%%%%%%%%%%%%%%%%%%%%%%%%%%%%%%%%%%%%%%%%%%%%%%%%%%%%%%%%%%%%%%%%%%%%%%%%%%%%%%%%%%%%%%%%%%%%%%%%%%%%%%%

Let $H$ be a graph.
We let $\C^{*}(H)$ denote the set of those components $C\in \C_{\rm odd}(H)$ such that $C$ is a hypomatchable graph having no $\{P_{2},P_{9}\}$-factor, and let $\C^{*}_{\leq 5}(H)=\{C\in \C^{*}(H)\mid |V(C)|\leq 5\}$, $\C^{*}_{\geq 7}(H)=\{C\in \C^{*}(H)\mid |V(C)|\geq 7\}$ and $\C^{**}_{\geq 7}(H)=\{C\in \C^{*}_{\geq 7}(H)\mid C$ is isomorphic to $K_{1}+sK_{2}$ for some $s\geq 3\}$.

\medbreak\noindent\textit{Proof of Theorem~\ref{mainthm2}.}\quad
Let $G$ be as in Theorem~\ref{mainthm2}.
Choose $S\subseteq V(G)$ so that (S1) and (S2) hold.

Set $T=\C_{\rm odd}(G-S)~(=\C(G-S))$, $T_{1}=\C^{*}_{\leq 5}(G-S)$ and $T_{2}=\C^{*}_{\geq 7}(G-S)$.
Then $T_{1}\cap T_{2}=\emptyset $ and $T_{1}\cup T_{2}\subseteq T$.
Now we construct a bipartite graph $H$ with bipartition $(S,T)$ by letting $uC\in E(H)~(u\in S,C\in T)$ if and only if $N_{G}(u)\cap V(C)\not=\emptyset $.
Let $L$ be the set of those edges $uC\in E(H)$ such that $u\in S$, $C\in \C^{**}_{\geq 7}(G-S)$ and $N_{G}(u)\cap V(C)$ consists only of the unique cutvertex of $C$.

\begin{claim}%%%%%%%%%%%%%%%%%%%%%%%%%%%%%%%%%%%%%%%%%%%%%%%%%%%%%%%%%%%%%%%%%%%%%%%%%%%%%%%%%%%%%%%%%%%%%%%%%%%%%%%%%%%
\label{cl5-2.1}
For every $Y\subseteq T_{1}\cup T_{2}$, $|N_{H-L}(Y)|\geq |Y\cap T_{1}|+\frac{1}{2}|Y\cap T_{2}|$.
\end{claim}
%%%%%%%%%%%%%%%%%%%%%%%%%%%%%%%%%%%%%%%%%%%%%%%%%%%%%%%%%%%%%%%%%%%%%%%%%%%%%%%%%%%%%%%%%%%%%%%%%%%%%%%%%%%%%%%%%%%%%%%%
\proof
Suppose that there exists $Y\subseteq T_{1}\cup T_{2}$ such that $|N_{H-L}(Y)|<|Y\cap T_{1}|+\frac{1}{2}|Y\cap T_{2}|$.
Set $X'=N_{H-L}(Y)$.
%Then each element of $Y\cap T_{1}$ belongs to $\C^{*}_{\leq 5}(G-X')$, and hence $|Y\cap T_{1}|\leq |\C^{*}_{\leq 5}(G-X')|\leq \sum _{0\leq j\leq 2}c_{2j+1}(G-X')$.
We divide $Y\cap T_{2}$ into two disjoint sets.
Let $Z_{1}$ be the set of those elements $C$ of $Y\cap T_{2}$ such that $|V(C)|=7$ and $C\not\in \C^{**}_{\geq 7}(G-S)$, and let $Z_{2}=(Y\cap T_{2})-Z_{1}$.
Note that $Z_{2}$ is the set of those elements $C$ of $Y\cap T_{2}$ such that $C$ is either isomorphic to $K_{1}+3K_{2}$ or a hypomatchable graph of order at least $9$ with no $\{P_{2},P_{9}\}$-factor.
Hence by the definition of $\G_{0}$ and Proposition~\ref{prop4-4.1}, each element of $Z_{2}$ belongs to $\G_{0}\cup \G_{1}\cup \G_{2}\cup \G_{3}\cup \G_{4}$.
Write $Z_{2}=\{D_{1},\ldots ,D_{q}\}$.
Let $X_{i}$ be a crush set of $D_{i}$ for each $1\leq i\leq q$, and set $X_{0}=X'\cup (\bigcup _{1\leq i\leq q}X_{i})$.
Let $1\leq i\leq q$.
We show that $\bigcup _{0\leq j\leq 2}\C_{2j+1}(D_{i}-X_{i})\subseteq \bigcup _{0\leq j\leq 2}\C_{2j+1}(G-X_{0})$.
This clearly holds if $D_{i}$ is a component of $G-X'$.
Thus we may assume that $D_{i}$ is not a component of $G-X'$.
By the definition of $L$, this means that $D_{i}\in \C^{**}_{\geq 7}(G-S)$ and the unique cutvertex of $D_{i}$ is the only vertex of $D_{i}$ that is adjacent to vertices in $S-X'$.
On the other hand, the unique cutvertex of $D_{i}$ is contained in $X_{i}$ by the definition of a crush set.
Hence $\bigcup _{0\leq j\leq 2}\C_{2j+1}(D_{i}-X_{i})\subseteq \bigcup _{0\leq j\leq 2}\C_{2j+1}(G-X_{0})$.

Since $i$ is arbitrary, we see that $c_{2j+1}(G-X_{0})=c_{2j+1}(G-X')+\sum _{1\leq i\leq q}c_{2j+1}(D_{i}-X_{i})$ for each $0\leq j\leq 2$.
By Lemma~\ref{lem4-2.1}, $c_{1}(D_{i}-X_{i})+c_{3}(D_{i}-X_{i})+\frac{2}{3}c_{5}(D_{i}-X_{i})\geq \frac{3}{4}|X_{i}|$ and $|X_{i}|\geq 4$ for every $1\leq i\leq q$.
Consequently
\begin{align*}
c_{1}(G-X_{0})&+c_{3}(G-X_{0})+\frac{2}{3}c_{5}(G-X_{0})\\
&\geq c_{1}(G-X')+c_{3}(G-X')+\frac{2}{3}c_{5}(G-X')+\frac{3}{4}\sum _{1\leq i\leq q}|X_{i}|.
\end{align*}
Hence
\begin{align}
c_{1}(G-X_{0})&+c_{3}(G-X_{0})+\frac{2}{3}c_{5}(G-X_{0})-\frac{2}{3}\sum _{0\leq j\leq 2}c_{2j+1}(G-X')-\frac{1}{3}q-\frac{2}{3}\sum _{1\leq i\leq q}|X_{i}|\nonumber \\
&\geq c_{1}(G-X')+c_{3}(G-X')+\frac{2}{3}c_{5}(G-X')+\frac{3}{4}\sum _{1\leq i\leq q}|X_{i}|\nonumber \\
&\quad -\frac{2}{3}\sum _{0\leq j\leq 2}c_{2j+1}(G-X')-\frac{1}{3}q-\frac{2}{3}\sum _{1\leq i\leq q}|X_{i}|\nonumber \\ 
&= \frac{1}{3}c_{1}(G-X_{0})+\frac{1}{3}c_{3}(G-X_{0})+\sum _{1\leq i\leq q}\left(\frac{3}{4}|X_{i}|-\frac{1}{3}-\frac{2}{3}|X_{i}|\right)\nonumber \\
&= \frac{1}{3}c_{1}(G-X_{0})+\frac{1}{3}c_{3}(G-X_{0})+\sum _{1\leq i\leq q}\left(\frac{1}{12}|X_{i}|-\frac{1}{3}\right)\nonumber \\
&\geq 0,\nonumber
\end{align}
which implies
$$
\frac{2}{3}\sum _{0\leq j\leq 2}c_{2j+1}(G-X')+\frac{1}{3}q+\frac{2}{3}\sum _{1\leq i\leq q}|X_{i}|\leq c_{1}(G-X_{0})+c_{3}(G-X_{0})+\frac{2}{3}c_{5}(G-X_{0}).
$$

Recall the definition of $X'$, $Z_{1}$ and $X_{0}$.
Since each element of $Y\cap T_{1}$ belongs to $\C^{*}_{\leq 5}(G-X')$, we have $|Y\cap T_{1}|\leq |\C^{*}_{\leq 5}(G-X')|\leq \sum _{0\leq j\leq 2}c_{2j+1}(G-X')$.
Since each element of $Z_{1}$ belongs to $\C_{7}(G-X_{0})$, we have $|Z_{1}|\leq c_{7}(G-X_{0})$.
Therefore
\begin{align*}
\frac{2}{3}|X_{0}| &= \frac{2}{3}\left(|X'|+\sum _{1\leq i\leq q}|X_{i}|\right)\\
&= \frac{2}{3}\left(|N_{H-L}(Y)|+\sum _{1\leq i\leq q}|X_{i}|\right)\\
&< \frac{2}{3}\left(|Y\cap T_{1}|+\frac{1}{2}|Y\cap T_{2}|+\sum _{1\leq i\leq q}|X_{i}|\right)\\
&\leq \frac{2}{3}\left(\sum _{0\leq j\leq 2}c_{2j+1}(G-X')+\frac{1}{2}(|Z_{1}|+|Z_{2}|)+\sum _{1\leq i\leq q}|X_{i}|\right)\\
&\leq \frac{2}{3}\left(\sum _{0\leq j\leq 2}c_{2j+1}(G-X')+\frac{1}{2}(c_{7}(G-X_{0})+q)+\sum _{1\leq i\leq q}|X_{i}|\right)\\
&\leq c_{1}(G-X_{0})+c_{3}(G-X_{0})+\frac{2}{3}c_{5}(G-X_{0})+\frac{1}{3}c_{7}(G-X_{0}),
\end{align*}
which contradicts the assumption of the theorem.
\qed

Now we apply Proposition~\ref{prop3.1} with $G$ replaced by $H$.
Then by Lemma~\ref{lem5.1}(iii) and Claim~\ref{cl5-2.1}, $H$ has a subgraph $F$ with $V(F)\supseteq S\cup T_{1}\cup T_{2}$ such that each $A\in \C(F)$ is a path satisfying one of (I) and (II) in Proposition~\ref{prop3.1}.
For $A\in \C(F)$, let $U_{A}=V(A)\cap S$ and $\L_{A}=V(A)\cap T$, and let $G_{A}=G[U_{A}\cup (\bigcup _{C\in \L_{A}}V(C))]$.

\begin{claim}%%%%%%%%%%%%%%%%%%%%%%%%%%%%%%%%%%%%%%%%%%%%%%%%%%%%%%%%%%%%%%%%%%%%%%%%%%%%%%%%%%%%%%%%%%%%%%%%%%%%%%%%%%%
\label{cl5-2.2}
For each $A\in \C(F)$, $G_{A}$ has a $\{P_{2},P_{9}\}$-factor.
\end{claim}
%%%%%%%%%%%%%%%%%%%%%%%%%%%%%%%%%%%%%%%%%%%%%%%%%%%%%%%%%%%%%%%%%%%%%%%%%%%%%%%%%%%%%%%%%%%%%%%%%%%%%%%%%%%%%%%%%%%%%%%%
\proof
We first assume that $A$ satisfies (I).
Then $|U_{A}|=|\L_{A}|=1$.
Write $U_{A}=\{u\}$ and $\L_{A}=\{D\}$, and let $v\in V(D)$ be a vertex with $uv\in E(G)$.
Since $D$ is hypomatchable by Lemma~\ref{lem5.1}(ii), $D-v$ has a perfect matching $M$.
Hence $M\cup \{uv\}$ is a perfect matching of $G_{A}$.
In particular, $G_{A}$ has a $\{P_{2},P_{9}\}$-factor.

Next we assume that $A$ satisfies (II).
Note that $|V(A)|$ is odd and $|V(A)|\geq 3$.
Write $A=D_{1}u_{1}D_{2}u_{2}\cdots D_{l}u_{l}D_{l+1}~(u_{i}\in U_{A},D_{i}\in \L_{A})$.
For $1\leq i\leq l$, let $v_{i}\in N_{G}(u_{i})\cap V(D_{i})$ and $w_{i+1}\in N_{G}(u_{i})\cap V(D_{i+1})$.
Since $u_{1}D_{1}$ and $u_{l}D_{l+1}$ are edges of $H-L$, we may assume that $v_{1}$ is not the unique cutvertex of $D_{1}$ if $D_{1}\simeq K_{1}+sK_{2}$ for some $s\geq 3$, and $w_{l+1}$ is not the unique cutvertex of $D_{l+1}$ if $D_{l+1}\simeq K_{1}+s'K_{2}$ for some $s'\geq 3$.
Since $D_{1}$ and $D_{l+1}$ are hypomatchable graphs of order at least $7$ by the definition of $T_{2}$, it follows from Lemma~\ref{lem4-5.2} that $D_{1}$ contains a path $Q_{1}$ with endvertex $v_{1}$ such that $|V(Q_{1})|\geq 5$ and $D_{1}-V(Q_{1})$ has a perfect matching $M_{1}$, and $D_{l+1}$ contains a path $Q_{l+1}$ with endvertex $w_{l+1}$ such that $|V(Q_{l+1})|\geq 5$ and $D_{l+1}-V(Q_{l+1})$ has a perfect matching $M_{l+1}$.
We regard $v_{1}$ as the terminal vertex of $Q_{1}$, and $w_{l+1}$ as the initial vertex of $Q_{l+1}$.
For each $i~(2\leq i\leq l)$, since $D_{i}$ is hypomatchable by the definition of $T_{1}$, it follows from Lemma~\ref{lem4-5.1} that $D_{i}$ contains a path $Q_{i}$ connecting $w_{i}$ to $v_{i}$ such that $D_{i}-V(Q_{i})$ has a perfect matching $M_{i}$.
Hence $P=Q_{1}u_{1}Q_{2}u_{2}\cdots Q_{l}u_{l}Q_{l+1}$ is a path of $G_{A}$ having order at least $11$.
Consequently $F_{A}=P\cup (\bigcup _{1\leq i\leq l+1}M_{i})$ is a path-factor of $G_{A}$ with $\C_{3}(F_{A})=\C_{5}(F_{A})=\C_{7}(F_{A})=\emptyset $ (and $\C_{9}(F_{A})=\emptyset $).
By Fact~\ref{fact1}, $G_{A}$ has a $\{P_{2},P_{9}\}$-factor.
\qed

By Lemma~\ref{lem5.1}(i)(ii), each component in $\C(G-S)-\C^{*}_{\leq 5}(G-S)-\C^{*}_{\geq 7}(G-S)$ has a $\{P_{2},P_{9}\}$-factor.
This together with Claim~\ref{cl5-2.2} implies that $G$ has a $\{P_{2},P_{9}\}$-factor.

This completes the proof of Theorem~\ref{mainthm2}.
\qed

%%%%%%%%%%%%%%%%%%%%%%%%%%%%%%%%%%%%%%%%%%%%%%%%%%%%%%%%%%%%%%%%%%%%%%%%%%%%%%%%%%%%%%%%%%%%%%%%%%%%%%%%%%%%%%%%%%%%%%%%
%%%%%%%%%%%%%%%%%%%%%%%%%%%%%%%%%%%%%%%%%%%%%%%%%%%%%%%%%%%%%%%%%%%%%%%%%%%%%%%%%%%%%%%%%%%%%%%%%%%%%%%%%%%%%%%%%%%%%%%%
%%%%%%%%%%%%%%%%%%%%%%%%%%%%%%%%%%%%%%%%%%%%%%%%%%%%%%%%%%%%%%%%%%%%%%%%%%%%%%%%%%%%%%%%%%%%%%%%%%%%%%%%%%%%%%%%%%%%%%%%
\section{Sharpness of Theorems~\ref{mainthm1} and \ref{mainthm2}}\label{sec6}
%%%%%%%%%%%%%%%%%%%%%%%%%%%%%%%%%%%%%%%%%%%%%%%%%%%%%%%%%%%%%%%%%%%%%%%%%%%%%%%%%%%%%%%%%%%%%%%%%%%%%%%%%%%%%%%%%%%%%%%%
%%%%%%%%%%%%%%%%%%%%%%%%%%%%%%%%%%%%%%%%%%%%%%%%%%%%%%%%%%%%%%%%%%%%%%%%%%%%%%%%%%%%%%%%%%%%%%%%%%%%%%%%%%%%%%%%%%%%%%%%
%%%%%%%%%%%%%%%%%%%%%%%%%%%%%%%%%%%%%%%%%%%%%%%%%%%%%%%%%%%%%%%%%%%%%%%%%%%%%%%%%%%%%%%%%%%%%%%%%%%%%%%%%%%%%%%%%%%%%%%%

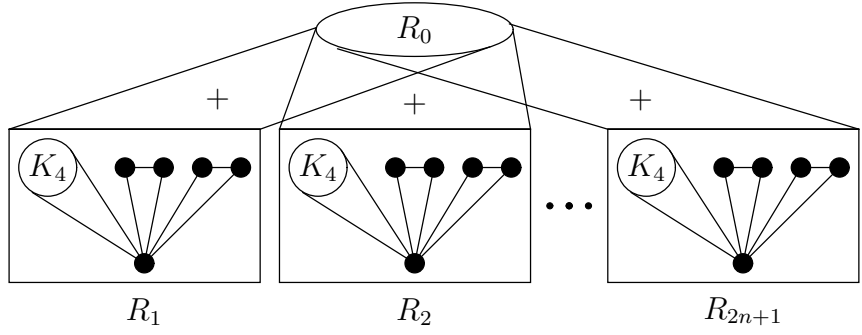
\begin{figure}
\begin{center}
%WinTpicVersion4.28b
{\unitlength 0.1in
\begin{picture}( 44.0000, 15.4500)(  2.0000,-18.8500)
% DOT 0 0 3 0 Black Black
% 4 3000 1400 3200 1400 3100 1400 3100 1400
% 
\special{pn 4}%
\special{sh 1}%
\special{ar 3000 1400 16 16 0  6.28318530717959E+0000}%
\special{sh 1}%
\special{ar 3200 1400 16 16 0  6.28318530717959E+0000}%
\special{sh 1}%
\special{ar 3100 1400 16 16 0  6.28318530717959E+0000}%
\special{sh 1}%
\special{ar 3100 1400 16 16 0  6.28318530717959E+0000}%
% CIRCLE 2 0 3 0 Black Black
% 4 400 1200 550 1200 550 1200 550 1200
% 
\special{pn 8}%
\special{ar 400 1200 150 150  0.0000000  6.2831853}%
% CIRCLE 2 0 0 0 Black Black
% 4 800 1200 800 1250 800 1250 800 1250
% 
\special{sh 1.000}%
\special{ia 800 1200 50 50  0.0000000  6.2831853}%
\special{pn 8}%
\special{ar 800 1200 50 50  0.0000000  6.2831853}%
% CIRCLE 2 0 0 0 Black Black
% 4 1000 1200 1000 1250 1000 1250 1000 1250
% 
\special{sh 1.000}%
\special{ia 1000 1200 50 50  0.0000000  6.2831853}%
\special{pn 8}%
\special{ar 1000 1200 50 50  0.0000000  6.2831853}%
% LINE 2 0 3 0 Black Black
% 2 1000 1200 800 1200
% 
\special{pn 8}%
\special{pa 1000 1200}%
\special{pa 800 1200}%
\special{fp}%
% CIRCLE 2 0 0 0 Black Black
% 4 1200 1200 1200 1250 1200 1250 1200 1250
% 
\special{sh 1.000}%
\special{ia 1200 1200 50 50  0.0000000  6.2831853}%
\special{pn 8}%
\special{ar 1200 1200 50 50  0.0000000  6.2831853}%
% CIRCLE 2 0 0 0 Black Black
% 4 1400 1200 1400 1250 1400 1250 1400 1250
% 
\special{sh 1.000}%
\special{ia 1400 1200 50 50  0.0000000  6.2831853}%
\special{pn 8}%
\special{ar 1400 1200 50 50  0.0000000  6.2831853}%
% LINE 2 0 3 0 Black Black
% 2 1400 1200 1200 1200
% 
\special{pn 8}%
\special{pa 1400 1200}%
\special{pa 1200 1200}%
\special{fp}%
% CIRCLE 2 0 0 0 Black Black
% 4 900 1700 900 1750 900 1750 900 1750
% 
\special{sh 1.000}%
\special{ia 900 1700 50 50  0.0000000  6.2831853}%
\special{pn 8}%
\special{ar 900 1700 50 50  0.0000000  6.2831853}%
% LINE 2 0 3 0 Black Black
% 4 800 1200 900 1700 900 1700 1000 1200
% 
\special{pn 8}%
\special{pa 800 1200}%
\special{pa 900 1700}%
\special{fp}%
\special{pa 900 1700}%
\special{pa 1000 1200}%
\special{fp}%
% LINE 2 0 3 0 Black Black
% 4 1200 1200 900 1700 900 1700 1400 1200
% 
\special{pn 8}%
\special{pa 1200 1200}%
\special{pa 900 1700}%
\special{fp}%
\special{pa 900 1700}%
\special{pa 1400 1200}%
\special{fp}%
% STR 2 0 3 0 Black Black
% 4 400 1100 400 1200 5 0 0 0
% $K_{4}$
\put(4.0000,-12.0000){\makebox(0,0){$K_{4}$}}%
% LINE 2 0 3 0 Black Black
% 2 900 1700 520 1110
% 
\special{pn 8}%
\special{pa 900 1700}%
\special{pa 520 1110}%
\special{fp}%
% LINE 2 0 3 0 Black Black
% 2 900 1700 320 1330
% 
\special{pn 8}%
\special{pa 900 1700}%
\special{pa 320 1330}%
\special{fp}%
% BOX 2 0 3 0 Black Black
% 2 200 1000 1500 1800
% 
\special{pn 8}%
\special{pa 200 1000}%
\special{pa 1500 1000}%
\special{pa 1500 1800}%
\special{pa 200 1800}%
\special{pa 200 1000}%
\special{pa 1500 1000}%
\special{fp}%
% CIRCLE 2 0 3 0 Black Black
% 4 1800 1200 1950 1200 1950 1200 1950 1200
% 
\special{pn 8}%
\special{ar 1800 1200 150 150  0.0000000  6.2831853}%
% CIRCLE 2 0 0 0 Black Black
% 4 2200 1200 2200 1250 2200 1250 2200 1250
% 
\special{sh 1.000}%
\special{ia 2200 1200 50 50  0.0000000  6.2831853}%
\special{pn 8}%
\special{ar 2200 1200 50 50  0.0000000  6.2831853}%
% CIRCLE 2 0 0 0 Black Black
% 4 2400 1200 2400 1250 2400 1250 2400 1250
% 
\special{sh 1.000}%
\special{ia 2400 1200 50 50  0.0000000  6.2831853}%
\special{pn 8}%
\special{ar 2400 1200 50 50  0.0000000  6.2831853}%
% LINE 2 0 3 0 Black Black
% 2 2400 1200 2200 1200
% 
\special{pn 8}%
\special{pa 2400 1200}%
\special{pa 2200 1200}%
\special{fp}%
% CIRCLE 2 0 0 0 Black Black
% 4 2600 1200 2600 1250 2600 1250 2600 1250
% 
\special{sh 1.000}%
\special{ia 2600 1200 50 50  0.0000000  6.2831853}%
\special{pn 8}%
\special{ar 2600 1200 50 50  0.0000000  6.2831853}%
% CIRCLE 2 0 0 0 Black Black
% 4 2800 1200 2800 1250 2800 1250 2800 1250
% 
\special{sh 1.000}%
\special{ia 2800 1200 50 50  0.0000000  6.2831853}%
\special{pn 8}%
\special{ar 2800 1200 50 50  0.0000000  6.2831853}%
% LINE 2 0 3 0 Black Black
% 2 2800 1200 2600 1200
% 
\special{pn 8}%
\special{pa 2800 1200}%
\special{pa 2600 1200}%
\special{fp}%
% CIRCLE 2 0 0 0 Black Black
% 4 2300 1700 2300 1750 2300 1750 2300 1750
% 
\special{sh 1.000}%
\special{ia 2300 1700 50 50  0.0000000  6.2831853}%
\special{pn 8}%
\special{ar 2300 1700 50 50  0.0000000  6.2831853}%
% LINE 2 0 3 0 Black Black
% 4 2200 1200 2300 1700 2300 1700 2400 1200
% 
\special{pn 8}%
\special{pa 2200 1200}%
\special{pa 2300 1700}%
\special{fp}%
\special{pa 2300 1700}%
\special{pa 2400 1200}%
\special{fp}%
% LINE 2 0 3 0 Black Black
% 4 2600 1200 2300 1700 2300 1700 2800 1200
% 
\special{pn 8}%
\special{pa 2600 1200}%
\special{pa 2300 1700}%
\special{fp}%
\special{pa 2300 1700}%
\special{pa 2800 1200}%
\special{fp}%
% STR 2 0 3 0 Black Black
% 4 1800 1100 1800 1200 5 0 0 0
% $K_{4}$
\put(18.0000,-12.0000){\makebox(0,0){$K_{4}$}}%
% LINE 2 0 3 0 Black Black
% 2 2300 1700 1920 1110
% 
\special{pn 8}%
\special{pa 2300 1700}%
\special{pa 1920 1110}%
\special{fp}%
% LINE 2 0 3 0 Black Black
% 2 2300 1700 1720 1330
% 
\special{pn 8}%
\special{pa 2300 1700}%
\special{pa 1720 1330}%
\special{fp}%
% BOX 2 0 3 0 Black Black
% 2 1600 1000 2900 1800
% 
\special{pn 8}%
\special{pa 1600 1000}%
\special{pa 2900 1000}%
\special{pa 2900 1800}%
\special{pa 1600 1800}%
\special{pa 1600 1000}%
\special{pa 2900 1000}%
\special{fp}%
% CIRCLE 2 0 3 0 Black Black
% 4 3500 1200 3650 1200 3650 1200 3650 1200
% 
\special{pn 8}%
\special{ar 3500 1200 150 150  0.0000000  6.2831853}%
% CIRCLE 2 0 0 0 Black Black
% 4 3900 1200 3900 1250 3900 1250 3900 1250
% 
\special{sh 1.000}%
\special{ia 3900 1200 50 50  0.0000000  6.2831853}%
\special{pn 8}%
\special{ar 3900 1200 50 50  0.0000000  6.2831853}%
% CIRCLE 2 0 0 0 Black Black
% 4 4100 1200 4100 1250 4100 1250 4100 1250
% 
\special{sh 1.000}%
\special{ia 4100 1200 50 50  0.0000000  6.2831853}%
\special{pn 8}%
\special{ar 4100 1200 50 50  0.0000000  6.2831853}%
% LINE 2 0 3 0 Black Black
% 2 4100 1200 3900 1200
% 
\special{pn 8}%
\special{pa 4100 1200}%
\special{pa 3900 1200}%
\special{fp}%
% CIRCLE 2 0 0 0 Black Black
% 4 4300 1200 4300 1250 4300 1250 4300 1250
% 
\special{sh 1.000}%
\special{ia 4300 1200 50 50  0.0000000  6.2831853}%
\special{pn 8}%
\special{ar 4300 1200 50 50  0.0000000  6.2831853}%
% CIRCLE 2 0 0 0 Black Black
% 4 4500 1200 4500 1250 4500 1250 4500 1250
% 
\special{sh 1.000}%
\special{ia 4500 1200 50 50  0.0000000  6.2831853}%
\special{pn 8}%
\special{ar 4500 1200 50 50  0.0000000  6.2831853}%
% LINE 2 0 3 0 Black Black
% 2 4500 1200 4300 1200
% 
\special{pn 8}%
\special{pa 4500 1200}%
\special{pa 4300 1200}%
\special{fp}%
% CIRCLE 2 0 0 0 Black Black
% 4 4000 1700 4000 1750 4000 1750 4000 1750
% 
\special{sh 1.000}%
\special{ia 4000 1700 50 50  0.0000000  6.2831853}%
\special{pn 8}%
\special{ar 4000 1700 50 50  0.0000000  6.2831853}%
% LINE 2 0 3 0 Black Black
% 4 3900 1200 4000 1700 4000 1700 4100 1200
% 
\special{pn 8}%
\special{pa 3900 1200}%
\special{pa 4000 1700}%
\special{fp}%
\special{pa 4000 1700}%
\special{pa 4100 1200}%
\special{fp}%
% LINE 2 0 3 0 Black Black
% 4 4300 1200 4000 1700 4000 1700 4500 1200
% 
\special{pn 8}%
\special{pa 4300 1200}%
\special{pa 4000 1700}%
\special{fp}%
\special{pa 4000 1700}%
\special{pa 4500 1200}%
\special{fp}%
% STR 2 0 3 0 Black Black
% 4 3500 1100 3500 1200 5 0 0 0
% $K_{4}$
\put(35.0000,-12.0000){\makebox(0,0){$K_{4}$}}%
% LINE 2 0 3 0 Black Black
% 2 4000 1700 3620 1110
% 
\special{pn 8}%
\special{pa 4000 1700}%
\special{pa 3620 1110}%
\special{fp}%
% LINE 2 0 3 0 Black Black
% 2 4000 1700 3420 1330
% 
\special{pn 8}%
\special{pa 4000 1700}%
\special{pa 3420 1330}%
\special{fp}%
% BOX 2 0 3 0 Black Black
% 2 3300 1000 4600 1800
% 
\special{pn 8}%
\special{pa 3300 1000}%
\special{pa 4600 1000}%
\special{pa 4600 1800}%
\special{pa 3300 1800}%
\special{pa 3300 1000}%
\special{pa 4600 1000}%
\special{fp}%
% STR 2 0 3 0 Black Black
% 4 2300 380 2300 480 5 0 0 0
% $R_{0}$
\put(23.0000,-4.8000){\makebox(0,0){$R_{0}$}}%
% STR 2 0 3 0 Black Black
% 4 1280 740 1280 840 5 0 0 0
% $+$
\put(12.8000,-8.4000){\makebox(0,0){$+$}}%
% STR 2 0 3 0 Black Black
% 4 2300 770 2300 870 5 0 0 0
% $+$
\put(23.0000,-8.7000){\makebox(0,0){$+$}}%
% STR 2 0 3 0 Black Black
% 4 3470 750 3470 850 5 0 0 0
% $+$
\put(34.7000,-8.5000){\makebox(0,0){$+$}}%
% ELLIPSE 2 0 3 0 Black Black
% 4 2300 480 2810 620 2810 620 2810 620
% 
\special{pn 8}%
\special{ar 2300 480 510 140  0.0000000  6.2831853}%
% LINE 2 0 3 0 Black Black
% 4 200 1000 1790 470 1500 1000 2690 570
% 
\special{pn 8}%
\special{pa 200 1000}%
\special{pa 1790 470}%
\special{fp}%
\special{pa 1500 1000}%
\special{pa 2690 570}%
\special{fp}%
% LINE 2 0 3 0 Black Black
% 2 1600 1000 1790 470
% 
\special{pn 8}%
\special{pa 1600 1000}%
\special{pa 1790 470}%
\special{fp}%
% LINE 2 0 3 0 Black Black
% 2 2900 1000 2810 490
% 
\special{pn 8}%
\special{pa 2900 1000}%
\special{pa 2810 490}%
\special{fp}%
% LINE 2 0 3 0 Black Black
% 2 3300 1000 1900 580
% 
\special{pn 8}%
\special{pa 3300 1000}%
\special{pa 1900 580}%
\special{fp}%
% LINE 2 0 3 0 Black Black
% 2 4600 1000 2790 450
% 
\special{pn 8}%
\special{pa 4600 1000}%
\special{pa 2790 450}%
\special{fp}%
% STR 2 0 3 0 Black Black
% 4 900 1850 900 1950 5 0 0 0
% $R_{1}$
\put(9.0000,-19.5000){\makebox(0,0){$R_{1}$}}%
% STR 2 0 3 0 Black Black
% 4 2300 1850 2300 1950 5 0 0 0
% $R_{2}$
\put(23.0000,-19.5000){\makebox(0,0){$R_{2}$}}%
% STR 2 0 3 0 Black Black
% 4 4000 1850 4000 1950 5 0 0 0
% $R_{2n+1}$
\put(40.0000,-19.5000){\makebox(0,0){$R_{2n+1}$}}%
\end{picture}}%
\caption{Graph $H_{n}$}
\label{f-s6-1}
\end{center}
\end{figure}

We first consider the coefficient of $|X|$ in Theorem~\ref{mainthm2}.
Let $n\geq 1$ be an integer.
Let $R_{0}$ be a complete graph of order $n$.
For each $i~(1\leq i\leq 2n+1)$, let $R_{i}$ be a graph isomorphic to $K_{1}+(K_{4}\cup 2K_{2})$.
Let $H_{n}=R_{0}+(\bigcup _{1\leq i\leq 2n+1}R_{i})$ (see Figure~\ref{f-s6-1}).

For $1\leq i\leq 2n+1$, since $|V(R_{i})|=9$ and $R_{i}$ does not contain a path of order $9$, $R_{i}$ has no $\{P_{2},P_{9}\}$-factor.
Suppose that $H_{n}$ has a $\{P_{2},P_{9}\}$-factor $F$.
Then for each $i~(1\leq i\leq 2n+1)$, $F$ contains an edge joining $V(R_{i})$ and $V(R_{0})$.
Since $2n+1>2|V(R_{0})|$, this implies that there exists $x\in V(R_{0})$ such that $d_{F}(x)\geq 3$, which is a contradiction.
Thus $H_{n}$ has no $\{P_{2},P_{9}\}$-factor.

\begin{lem}%%%%%%%%%%%%%%%%%%%%%%%%%%%%%%%%%%%%%%%%%%%%%%%%%%%%%%%%%%%%%%%%%%%%%%%%%%%%%%%%%%%%%%%%%%%%%%%%%%%%%%%%%%%%%
\label{lem6.1}
For all $X\subseteq V(H_{n})$, $\sum _{0\leq j\leq 3}c_{2j+1}(H_{n}-X)\leq \frac{2}{3}|X|+\frac{1}{3}$.
\end{lem}
%%%%%%%%%%%%%%%%%%%%%%%%%%%%%%%%%%%%%%%%%%%%%%%%%%%%%%%%%%%%%%%%%%%%%%%%%%%%%%%%%%%%%%%%%%%%%%%%%%%%%%%%%%%%%%%%%%%%%%%%
\proof
Let $X\subseteq V(H_{n})$.

\begin{claim}%%%%%%%%%%%%%%%%%%%%%%%%%%%%%%%%%%%%%%%%%%%%%%%%%%%%%%%%%%%%%%%%%%%%%%%%%%%%%%%%%%%%%%%%%%%%%%%%%%%%%%%%%%%
\label{cl6.1.1}
For each $i~(1\leq i\leq 2n+1)$, $\sum _{0\leq j\leq 3}c_{2j+1}(R_{i}-X)\leq \frac{2}{3}|V(R_{i})\cap X|+\frac{1}{3}$.
\end{claim}
%%%%%%%%%%%%%%%%%%%%%%%%%%%%%%%%%%%%%%%%%%%%%%%%%%%%%%%%%%%%%%%%%%%%%%%%%%%%%%%%%%%%%%%%%%%%%%%%%%%%%%%%%%%%%%%%%%%%%%%%
\proof
Let $u$ be the unique cutvertex of $R_{i}$.

We first assume that $u\not\in X$.
Then $R_{i}-X$ is connected.
Clearly we may assume that $\sum _{0\leq j\leq 3}c_{2j+1}(R_{i}-X)=1$.
Then $|V(R_{i})\cap X|\geq 2$ because $|V(R_{i})|=9$.
Hence $\sum _{0\leq j\leq 3}c_{2j+1}(R_{i}-X)=1<\frac{2}{3}\cdot 2+\frac{1}{3}\leq \frac{2}{3}|V(R_{i})\cap X|+\frac{1}{3}$.
Thus we may assume that $u\in X$.

Let $\alpha $ be the number of components of $R_{i}-u$ intersecting with $X$.
Since $\alpha \leq 3$, we have $\alpha \leq \frac{2}{3}(\alpha +1)+\frac{1}{3}$.
Furthermore, $\sum _{0\leq j\leq 3}c_{2j+1}(R_{i}-X)=c_{1}(R_{i}-X)+c_{3}(R_{i}-X)\leq \alpha $ and $|V(R_{i})\cap X|=|\{u\}|+|(V(R_{i})-\{u\})\cap X|\geq \alpha +1$.
Consequently we get $\sum _{0\leq j\leq 3}c_{2j+1}(R_{i}-X)\leq \frac{2}{3}|V(R_{i})\cap X|+\frac{1}{3}$.
\qed

Assume for the moment that $V(R_{0})\not\subseteq X$.
Then $H_{n}-X$ is connected.
Clearly we may assume that $\sum _{0\leq j\leq 3}c_{2j+1}(H_{n}-X)=1$.
Then $|X|\geq 2$ because $|V(H_{n})|\geq 9$.
Hence $\sum _{0\leq j\leq 3}c_{2j+1}(H_{n}-X)=1<\frac{2}{3}\cdot 2+\frac{1}{3}\leq \frac{2}{3}|X|+\frac{1}{3}$.
Thus we may assume that $V(R_{0})\subseteq X$.
Then clearly
\begin{align}
|&\C_{2j+1}(H_{n}-X)|=\sum _{1\leq i\leq 2n+1}|\C_{2j+1}(R_{i}-X)|.\label{eq-sec6-1}
\end{align}
By Claim~\ref{cl6.1.1} and (\ref{eq-sec6-1}),
\begin{align*}
\sum _{0\leq j\leq 3}c_{2j+1}(H_{n}-X) &= \sum _{0\leq j\leq 3}\left(\sum _{1\leq i\leq 2n+1}c_{2j+1}(R_{i}-X)\right)\\
&\leq \sum _{1\leq i\leq 2n+1}\left(\frac{2}{3}|V(R_{i})\cap X|+\frac{1}{3}\right)\\
&= \frac{2}{3}(|X|-|V(R_{0})|)+\frac{1}{3}(2n+1)\\
&= \frac{2}{3}(|X|-n)+\frac{1}{3}(2n+1)\\
&= \frac{2}{3}|X|+\frac{1}{3}.
\end{align*}
Thus we get the desired conclusion.
\qed

From Lemma~\ref{lem6.1}, we get the following proposition, which implies that the coefficient of $|X|$ in Theorem~\ref{mainthm2} is best possible in the sense that it cannot be replaced by any number greater than $\frac{2}{3}$.

\begin{prop}%%%%%%%%%%%%%%%%%%%%%%%%%%%%%%%%%%%%%%%%%%%%%%%%%%%%%%%%%%%%%%%%%%%%%%%%%%%%%%%%%%%%%%%%%%%%%%%%%%%%%%%%%%%%
\label{prop3.2}
There exist infinitely many graphs $G$ having no $\{P_{2},P_{9}\}$-factor such that $\sum _{0\leq i\leq 3}c_{2i+1}(G-X)\leq \frac{2}{3}|X|+\frac{1}{3}$ for all $X\subseteq V(G)$.
\end{prop}
%%%%%%%%%%%%%%%%%%%%%%%%%%%%%%%%%%%%%%%%%%%%%%%%%%%%%%%%%%%%%%%%%%%%%%%%%%%%%%%%%%%%%%%%%%%%%%%%%%%%%%%%%%%%%%%%%%%%%%%%

We now briefly discuss the sharpness of other coefficients.
Let $n\geq 8$, and let $R_{0}$ be a complete graph of order $n$.
For each $i~(1\leq i\leq n+1)$, let $R_{i}$ be a graph isomorphic to $K_{1}+2K_{2}$, and let $u_{i}$ be the unique cutvertex of $R_{i}$.
Let $H$ be the graph obtained from $R_{0}\cup (\bigcup _{1\leq i\leq n+1}R_{i})$ by joining $u_{i}$ to all vertices in $R_{0}$ for each $i~(1\leq i\leq n+1)$.
Then $c_{1}(H-V(R_{0}))+c_{3}(H-V(R_{0}))+\frac{2}{3}c_{5}(H-V(R_{0}))+\frac{1}{3}c_{7}(H-V(R_{0}))=\frac{2}{3}c_{5}(H-V(R_{0}))=\frac{2}{3}|V(R_{0})|+\frac{2}{3}$, and $c_{1}(H-X)+c_{3}(H-X)+\frac{2}{3}c_{5}(H-X)+\frac{1}{3}c_{7}(H-X)\leq \frac{2}{3}|X|$ for all $X\subseteq V(H)$ with $X\not=V(R_{0})$, and $H$ has no $\{P_{2},P_{9}\}$-factor.
This shows that the coefficient of $c_{5}(G-X)$ in Theorem~\ref{mainthm2} is best possible in the sense that it cannot be replaced by any number less than $\frac{2}{3}$.
Similarly graphs $K_{n}+(2n+1)K_{7}~(n\geq 1)$ show that the coefficient of $c_{7}(G-X)$ in Theorem~\ref{mainthm2} is best possible in the sense that it cannot be replaced by any number less than $\frac{1}{3}$.

As for Theorem~\ref{mainthm1}, graphs $K_{n}+(2n+1)(K_{1}+3K_{2})~(n\geq 1)$ show that the coefficient $\frac{2}{3}$ of $|X|$ is best possible, and graphs $K_{n}+(2n+1)K_{3}$ and $K_{n}+(2n+1)K_{5}~(n\geq 1)$ show that the coefficient $\frac{1}{3}$ of $c_{3}(G-X)$ and $c_{5}(G-X)$ are best possible.

\end{document}